\documentclass{amsart}

\usepackage{xypic}
\input xy
\xyoption{all}
\usepackage{epsfig}
\usepackage{amsthm}
\usepackage{amssymb}
\usepackage{amsmath}
\usepackage{amscd}
\usepackage{color}

\newcommand{\bg}{\begin{equation}}
\newcommand{\ed}{\end{equation}}
\newcommand{\bga}{\begin{eqnarray}}
\newcommand{\eda}{\end{eqnarray}}
\newcommand{\pf}{\textbf{Proof:\ }}

\def\cbdu{\par{\raggedleft$\Box$\par}}

\newtheorem {Theorem}  {Theorem}

\numberwithin{Theorem}{section}

\newtheorem {Lemma}[Theorem]  {Lemma}

\theoremstyle{definition}

\theoremstyle{remark}
\newtheorem{Remark}[Theorem]{\bf Remark}

\expandafter\chardef\csname pre amssym.def
at\endcsname=\the\catcode`\@ \catcode`\@=11
\def\undefine#1{\let#1\undefined}
\def\newsymbol#1#2#3#4#5{\let\next@\relax
 \ifnum#2=\@ne\let\next@\msafam@\else
 \ifnum#2=\tw@\let\next@\msbfam@\fi\fi
 \mathchardef#1="#3\next@#4#5}
\def\mathhexbox@#1#2#3{\relax
 \ifmmode\mathpalette{}{\m@th\mathchar"#1#2#3}%
 \else\leavevmode\hbox{$\m@th\mathchar"#1#2#3$}\fi}
\def\hexnumber@#1{\ifcase#1 0\or 1\or 2\or 3\or 4\or 5\or 6\or 7\or 8\or
 9\or A\or B\or C\or D\or E\or F\fi}

\font\teneufm=eufm10 \font\seveneufm=eufm7 \font\fiveeufm=eufm5
\newfam\eufmfam
\textfont\eufmfam=\teneufm \scriptfont\eufmfam=\seveneufm
\scriptscriptfont\eufmfam=\fiveeufm

\catcode`\@=\csname pre amssym.def at\endcsname

\newcounter{remark}
\newcommand{\R}{\mathbf{R}}

\def  \R   {{\mathbb R}}

\def  \T   {{\mathbb T}}

\def  \12  {{\frac{1}{2}}}

\def\build#1_#2^#3{\mathrel{\mathop{\kern 0pt#1}\limits_{#2}^{#3}}}

\begin{document}

\title[Determining wavenumbers]{Determining wavenumbers for the  incompressible Hall-magneto-hydrodynamics}


\author [Han Liu]{Han Liu}
\address{Department of Mathematics, Stat. and Comp. Sci.,  University of Illinois Chicago, Chicago, IL 60607,USA}
\email{hliu94@uic.edu} 

\thanks{The work of the authors was partially supported by NSF Grant
DMS--1108864.}





\begin{abstract}
Using Littlewood-Paley theory, one formulates the determining wavenumbers for the Hall-MHD system, defined for each individual solution $(u,b)$. It is shown that the long time behaviour of strong solutions is almost finite dimensional as the wavenumbers are bounded in certain average senses.

\bigskip

KEY WORDS: Hall-MHD system; determining modes

\hspace{0.02cm} CLASSIFICATION CODE: 35Q35, 35Q85, 37L30
\end{abstract}

\maketitle

\section{Introduction}

This paper deals with the finite dimensionality of solutions to the  incompressible Hall-magneto-hydrodynamics (Hall-MHD) system, written as follows
\begin{align}
u_t+(u\cdot \nabla)u-(b \cdot \nabla)b+\nabla p =\nu \Delta u+f, \label{hmhd1}\\
b_t+(u \cdot \nabla )b-(b \cdot \nabla) u+\eta \nabla \times ((\nabla \times b) \times b)=\mu \Delta b, \label{hmhd2} \\
\nabla \cdot u=0,\ \nabla \cdot b=0, \ t \in \R^+, \ x \in \T^3. \label{hmhd3}
\end{align}
The above system describes the evolution of a system consisting of a magnetic field $b,$ electrons and ions, whose collective motion under $b$ can be approximated as an electrically conducting fluid with velocity field $u.$ The focus here is primarily the visco-resistive case, corresponding to positive fluid viscosity $\nu$ and magnetic resistivity $\mu.$ The external forcing term $f$, which one assumes to have zero mean, renders system (\ref{hmhd1})-(\ref{hmhd3}) inhomogeneous. The Hall-MHD system is derived using generalized Ohm's law which takes into account the effect of the electric current on the Lorentz force, neglected in the derivation of the MHD equations. The resulting extra Hall term $\eta(\nabla \times ((\nabla \times b) \times b)),$ distinguishing system (\ref{hmhd1})-(\ref{hmhd3}) from the conventional MHD system, becomes significant in the case of large magnetic shear. The coefficient $\eta$ here is proportional to the ion skin length.

The Hall-MHD system has a wide range of applications including modelling solar winds, designing magnetic confinement devices for fusion reactors and interpreting the origin of the geomagnetic field. It is believed to be an essential model for magnetic reconnection, an intriguing phenomenon frequently observed in space plasmas. Over the past decade, the Hall-MHD system has received more attentions from the mathematical community. Acheritogaray, Degond, Frouvelle and Liu \cite{ADFL} rigorously derived the system and established global existence of weak solutions on periodic domains. Chae, Degond and Liu \cite{CDL} proved existence of global weak solutions in $\R^3$ as well as that of local smooth solutions. Chae and Lee \cite{CL} obtained blow-up criteria and small data global strong solutions. Evidences of ill-posedness can be found in \cite{CW, D4, JO}. As for the properties of solutions, temporal decay estimates in energy spaces are due to Chae and Schonbek \cite{CS}. For more mathematical results on the Hall-MHD system, e.g., well-posedness results and regularity criteria, please see \cite{BF, D1, D2, D3, FFNZ, FLN, HAHZ, KL, WZ1, WL, Y3, Y4, YZ0, Z}.

In the case of $b\equiv 0,$ system (\ref{hmhd1})-(\ref{hmhd3}) reduces to the Navier-Stokes equations (NSE), for which the finite dimensional behaviour of solutions has been extensively studied. As alluded in Kolmogorov's 1941 phenomenological theory \cite{K}, a turbulent flow should have a finite number of degrees of freedom. The first mathematical result in this direction, due to Foia\c s and Prodi \cite{FP}, stated that the higher Fourier modes of a solution to the 2D NSE are controlled by the lower modes asymptotically as time goes to infinity. More precisely, if a certain finite number of Fourier modes of a solution share the same long time behaviour with those of another solution, then the remaining infinitely many Fourier modes of the two solutions also exhibit the same long time behaviour. Thus, the notion of ``determining modes" arises naturally. For the 2D NSE, estimates of the number of the determining modes were obtained by Foia\c s, Manley, Temam and Treve \cite{FMTT} in terms of the Grashof number, and later improved by Jones and Titi \cite{JT}, whereas Constantin, Foia\c s, Manley and Temam \cite{CFMT} estimated the number of determining modes for the 3D NSE assuming the uniform boundedness of solutions in $H^1$. For more details concerning the study of finite dimensionality of the NSE flow, readers are referred to \cite{CFT, FJKT, FMRT, FT1, FT2, FT3, R}.

Motivated by the work of Cheskidov, Dai and Kavlie \cite{CDK} where a time-dependent determining wavenumber was introduced to estimate the number of determining modes for weak solutions to the 3D NSE in an average sense, this paper aims to adapt the idea therein to the study of the Hall-MHD system. In particular, as finite dimensionality of the closely related MHD system has been investigated by Eden and Libin in \cite{EL}, one is curious if such results can also be obtained for the Hall-MHD system, which differs nontrivially from the MHD system in many aspects, as illustrated in \cite{CL, CW, JO}.

One introduces the determining wavenumbers for an individual weak solution $(u,b)$ to system (\ref{hmhd1})-(\ref{hmhd3}). Let $\kappa:= \min\{\mu, \nu, \eta^{-1}\mu\},$ $r \in (2,3)$ and $\delta >1.$ Let $c_r$ be a constant depending only on $r.$ The determining wavenumbers corresponding to $u$ and $b$ are defined as follows.
\begin{align}
\Lambda_{u}(t)=: & \min\big\{\lambda_q: \lambda_p^{-1+\frac{3}{r}}\|u_p\|_r < c_r\kappa, \forall p >q; \ \lambda_q^{-1+\frac{3}{r}}\|u_{\leq q}\|_r < c_r \kappa, q \in \mathbb{N} \big\}, \label{uvnbr}\\
\Lambda_{b}(t)=: & \min\big\{\lambda_q: \lambda^\delta_{p-q}\|b_p\|_\infty < c_r \kappa , \forall p >q; \ \|b_{\leq q}\|_\infty < c_r  \kappa, q \in \mathbb{N} \big\},\label{wvnbr}
\end{align}
where $\lambda_q = 2^q$ and $u_q = \Delta_q u,$ the $q$-th Littlewood-Paley projection of $u,$ to be defined in Section \ref{sec:pre}. One notices that in both definitions the conditions on the high modes resemble the ones in the definitions of the dissipation wavenumbers found in \cite{CD, D1}. It is noteworthy that the dissipation wavenumber, which separates the dissipation range from the inertial range of turbulent flows, has been utilized to establish improved regularity criteria for various fluid models in \cite{CD, CS, D1, D5}. In this paper, the following theorem shall be proved.
\begin{Theorem}\label{thm}
Let $(u,b)$ and $(v,h)$ be two weak solutions to system (\ref{hmhd1})-(\ref{hmhd3}) such that for all $t \geq 0,$ $$\frac{1}{|\T^3|}\int_{\T^3} u(t, x)\mathrm{d}x = \frac{1}{|\T^3|}\int_{\T^3} v(t, x) \mathrm{d}x =0 \text{ and }\frac{1}{|\T^3|}\int_{\T^3} b(t, x)-h(t, x)\mathrm{d}x =0.$$ Let $\Lambda_{u,v}(t):=\max\{\Lambda_{u}(t), \Lambda_{v}(t)\}$ and $\Lambda_{b,h}(t):=\max\{\Lambda_{b}(t), \Lambda_{h}(t)\}$ with $\Lambda_{u}(t),$ $\Lambda_{v}(t),$ $\Lambda_{b}(t)$ and $\Lambda_{h}(t)$ defined as in Definition \ref{uvnbr}-\ref{wvnbr}. Let $Q_{u,v}(t)$ and $Q_{b,h}(t)$ be such that $\Lambda_{u,v}(t)=\lambda_{Q_{u,v}(t)}\text{ and }\Lambda_{b,h}(t)=\lambda_{Q_{b,h}(t)}.$  If $$\big(u_{\leq Q_{u,v}(t)}(t), \ b_{\leq Q_{b,h}(t)}(t)\big)= \big(v_{\leq Q_{u,v}(t)}(t), h_{\leq Q_{b,h}(t)}(t)\big), \ \forall t>0,$$ then $$\lim_{t \to \infty} \big(\|u(t)-v(t)\|_{L^2}+\|b(t)-h(t)\|_{L^2}\big) =0.$$
\end{Theorem}

\begin{Remark} 
Due to the Galilean invariance of the fluid equation, it suffices to assume that $u$ and $v$ are zero-mean solutions. Yet, in general the magnetic fields do not have zero means, so in the above theorem it is assumed that the space-average of $b$ is the same as that of $h.$
\end{Remark}

\bigskip

\section{Preliminaries}
\label{sec:pre}

\subsection{Notation}
\label{sec:notation}

Throughout the paper, $A\lesssim B$ denotes an estimate of the form $A\leq C B$ with some absolute constant $C$. Given a tempered distribution $u,$ one denotes by $\mathcal{F}u=\hat u$ and $\mathcal{F}^{-1}u= \check u$ the Fourier transform and the inverse Fourier transform of $u,$ respectively. For simplicity, the $L^p$-norm $\|\cdot\|_{L^p}$ is sometimes written as $\|\cdot\|_p$ , while $H^{s}$ denotes the $L^2$-based Sobolev spaces.

\subsection{Well-posedness results for system (\ref{hmhd1})-(\ref{hmhd3})}
In order to discuss the determining modes of the solutions, one had better first clarify the notions and existence of the solutions. Relevant to this paper are the following well-posedness results. From \cite{ADFL}, it is known that in $\T^3,$ the Leray-Hopf type weak solution to system (\ref{hmhd1})-(\ref{hmhd3}) exists globally in time, just as in the case of the NSE.
\begin{Theorem}[Leray-Hopf type weak solution]
Let the initial data $u_0, b_0 \in L^2(\T^3).$ There exists a global weak solution $(u,b)$ to system (\ref{hmhd1})-(\ref{hmhd3}) satisfying 
$$(u,b) \in L^\infty \big(0,T; (L^2(\T^3))^2\big)\cap L^2\big(0,T; (H^1(\T^3))^2\big).$$ In addition, the following energy inequality holds -
$$\frac{1}{2}\frac{\mathrm{d}}{\mathrm{d}t}\big(\|u\|_{L^2}^2 + \|b\|_{L^2}^2\big) \leq - \nu \|\nabla u\|_{L^2}^2 -\mu \|\nabla b\|_{L^2}^2.$$ 
\end{Theorem}

In \cite{CDL}, local existence of strong solutions was proven. Furthermore, for small initial data, the existence of strong solutions is global.
\begin{Theorem}[Strong solution]
Let $s >\frac{5}{2}$ be an integer and $u_0, b_0 \in H^{s}(\T^3)$ with $\nabla \cdot u_0 = \nabla \cdot b_0 =0.$ Then:

i) The initial value $(u_0, b_0)$ generates a local-in-time classical solution $(u,b) \in L^\infty\big(0,T; (H^s(\T^3)^2)\big)$ with $T=T(\|u_0\|_{H^s}, \|b_0\|_{H^s}).$

ii) There exists a constant $\epsilon=\epsilon(\nu, s)$ such that $(u_0, b_0)$ generates a global classical solution $(u,b) \in L^\infty\big(0,\infty; (H^s(\T^3)^2)\big)$, provided that $\| u_0\|_{H^s}+\|b_0\|_{H^s}<\epsilon.$
\end{Theorem}

To demonstrate the finite dimensional behaviour of the solutions, the following regularity criterion, found in \cite{CL}, is needed.
\begin{Theorem}[Prodi-Serrin type regularity criterion]\label{rgc}
Let $s >\frac{5}{2}$ be an integer and $u_0, b_0 \in H^{s}(\T^3)$ with $\nabla \cdot u_0 = \nabla \cdot b_0 =0.$ Then for the first blow-up time $T^* <\infty$ of the classical solution to system (\ref{hmhd1})- (\ref{hmhd3}), it holds that 
$$\limsup_{t \nearrow T^*} (\|u(t)\|_{H^s}^2 + \|b(t)\|_{H^s}^2)=\infty,$$
if and only if $$\|u\|_{L^q(0,T^*; L^p(\T^3))}+\|\nabla b\|_{L^\gamma(0,T^*; L^\beta(\T^3))}=\infty,$$ 
where $p,q, \beta$ and $\gamma$ satisfy the relation
$$\frac{3}{p}+\frac{2}{q} \leq 1, \  \frac{3}{\beta}+\frac{2}{\gamma} \leq 1, \text{ with } p, \beta \in (3, \infty].$$
\end{Theorem}

\subsection{Littlewood-Paley decomposition}

This section is a brief introduction to the Littlewood-Paley theory, a fundamental tool used throughout the paper. One starts with introducing a family of functions with annular support, $\{\varphi_q(\xi)\}_{q=-1}^\infty$, which forms a dyadic partition of unity in the frequency domain. Let $\lambda_q =2^q, q \in \mathbb{Z}.$ One chooses a radial function $\chi \in C^\infty_0(\R^n)$ satisfying
\begin{equation}\notag \chi(\xi)=
\begin{cases}
1, \ \text{for }|\xi| \leq \frac{3}{4}\\
0, \ \text{for }|\xi| \geq 1,\\
\end{cases}
\end{equation}
and define $\varphi(x)=\chi(\frac{\xi}{2})-\chi(\xi)$ and $\varphi_q(\xi)=\begin{cases} \varphi(\lambda_q^{-1}\xi), \text{for }q \geq 0,  \\ \chi(\xi), \ \text{for }q= -1. \end{cases}$ 

For a vector field $u \in \mathcal{S}^{'}(\T^n),$ one defines the Littlewood-Paley projections as 
\begin{equation}\notag
\Delta_q u = u_q= : \sum_{k\in \mathbb{Z}^n} \varphi_q(k)\hat u(k)e^{i2\pi k\cdot x},
\end{equation}
where $\hat u(k)$ is the $k$-th Fourier coefficient of $u.$ In particular, $\hat u(0)=u_{-1}.$ Thus, at least in the distributional sense $u$ can be identified as a sum of its Littlewood-Paley projections $$u =\sum_{q=-1}^\infty u_q.$$ One also introduces the following notations, which appear throughout the paper, $$u_{ \leq Q}:= \sum_{q =-1}^Q u_q, \ \ u_{(P, Q]}:= \sum_{q= P+1}^Q u_q, \ \  \tilde u_q:= \sum_{|p-q| \leq 1} u_p.$$ 

The $L^2$-based Sobolev spaces can thus be characterized via  Littlewood-Paley projections -
$$\|u\|_{H^s} = \Big(\sum_{q \geq -1} \lambda_q^{2s}\|u_q\|_2^2\Big)^{\frac{1}{2}}.$$

In addition, the following Bernstein's inequality shall be used extensively.
\begin{Lemma}\label{brn}Let $n$ be the space dimension and let $s \geq r \geq 1,$ then
\begin{equation}\notag
\|u_q\|_{r} \lesssim \lambda_q^{n(\frac{1}{r}-\frac{1}{s})} \|u_q\|_s. 
\end{equation}
\end{Lemma}
\pf See \cite{BCD}.
\cbdu

\subsection{Bony's paraproduct and commutator estimates}

The product of two distributions $u$ and $v$ can be formally written as $$uv=\sum_{p,q \geq -1} u_p v_q.$$ Using Bony's paradifferential calculus, one has the following paraproduct decomposition
$$uv= \sum_{q \geq -1} u_{\leq q-2} v_q+ \sum_{q \geq -1} u_q v_{\leq q-2}+\sum_{q \geq -1} \tilde u_q v_q,$$ which distinguishes three parts in the product $uv.$

To facilitate the estimations, one introduces the following commutators for the convection or inertial terms and the Hall term, respectively
\begin{equation}\label{cm1}
[\Delta_q, u_{\leq p-2}\cdot \nabla] v_p = \Delta_q(u_{\leq p-2}\cdot \nabla v_p)-u_{\leq p-2} \cdot \nabla \Delta_q v_p,
\end{equation}
\begin{equation}\label{cm2}
[\Delta_q, b_{\leq p-2}\times \nabla \times] h_p = \Delta_q(b_{\leq p-2}\times (\nabla \times h_p))-b_{\leq p-2} \times (\nabla \times \Delta_q h_p).
\end{equation}

In the upcoming sections, it shall be seen that the commutators, along with the divergence free conditions, reveal certain cancellations within the nonlinear interactions. The commutator (\ref{cm1}) enjoys the following estimate, proven in \cite{BCD}.
\begin{Lemma}\label{cmest} For $\frac{1}{r_1}=\frac{1}{r_2}+\frac{1}{r_3},$ the following inequality is true -
\begin{equation}\notag
\|[\Delta_q, u_{\leq p-2}\cdot \nabla] v_p\|_{r_1} \lesssim \|v_p\|_{r_2} \sum_{p' \leq p-2}\lambda_{p'}\| u_{p'}\|_{r_3}.
\end{equation}
\end{Lemma}

The commutator (\ref{cm2}) satisfies an analogous estimate, as shown in \cite{D1}.
\begin{Lemma}\label{cmmes} Given that $\nabla \cdot b_{\leq p-2} =0,$ the following inequality holds -
\begin{equation}\notag
\|[\Delta_q, b_{\leq p-2}\times\nabla \times] h_p\|_{r} \lesssim \|h_p\|_{r} \sum_{p' \leq p-2}\lambda_{p'}\| b_{p'}\|_{\infty}.
\end{equation}
\end{Lemma}

More detailed study of the aforementioned harmonic analysis tools and their applications can be found in the work of Bahouri, Chemin and Danchin \cite{BCD}. 

\bigskip

\section{Analysis of a reduced system}\label{EMH}
To analyze the complete Hall-MHD system, one starts by considering the fluid-free version of system (\ref{hmhd1})-(\ref{hmhd3}), written as follows - 
\begin{align}
b_t + \eta \nabla \times ((\nabla \times b) \times b)=\mu \Delta b, \label{emhd1}\\
\nabla \cdot b =0. \label{emhd2}
\end{align}

The above system is named electron magneto-hydrodynamic (EMHD) equations as it describes the situation where the ions in the Hall-MHD setting are too heavy to move, leaving only the electrons in motion. As the small-scale limit of the Hall-MHD system, the EMHD equations can be used as a toy model to better understand the Hall term. In \cite{ADFL}, the existence of weak solutions to (\ref{emhd1})-(\ref{emhd2}) on periodic domains, analogous to that of the complete Hall-MHD system, was shown. For more studies concerning the EMHD equations, readers may consult \cite{G1, G3, MG}. 

In the following passages, $b$ and $h$ are two weak solutions to system (\ref{emhd1})-(\ref{emhd2}). One aims to prove the following analogue of Theorem \ref{thm}. 
\begin{Theorem}\label{trm}
Let $\Lambda_{b,h}(t):=\max\{\Lambda_{b}(t), \Lambda_{h}(t)\}$ with $\Lambda_{b}(t)$ and $\Lambda_{h}(t)$ defined as in definition (\ref{wvnbr}). Let $Q_{b,h}(t)$ be such that $\Lambda_{b,h}(t)=\lambda_{Q_{b,h}(t)}.$  If $$b_{\leq Q_{b,h}(t)}(t)= h_{\leq Q_{b,h}(t)}(t), \ \forall t>0,$$ then $$\lim_{t \to \infty} \|b(t)-h(t)\|_{L^2} =0, \ \forall s>0.$$
\end{Theorem} 

\pf Straightforward calculations show that $m:=b-h$ satisfies 
\begin{equation}\notag
m_t - \mu \Delta m = - \eta \nabla \times ((\nabla \times m) \times h) - \eta \nabla \times ((\nabla \times b) \times m).
\end{equation}

Multiplying the above equation by $\lambda_q^{2s}\Delta_q^2 m,$ integrating by parts and summing over $q$ lead to the following identity.
\begin{equation}\notag
\begin{split}
& \frac{1}{2}\frac{\mathrm{d}}{\mathrm{d}t}\sum_{q \geq -1}\lambda_q^{2s}\|m_q\|_2^2+\mu \sum_{q \geq -1}\lambda_q^{2s+2}\|m_q\|_2^2\\
= &\eta\sum_{q \geq -1} \lambda_q^{2s}\int_{\T^3} \Delta_q ( (\nabla \times m)\times h)\cdot (\nabla \times m_q) \mathrm{d}x\\
&+ \eta \sum_{q \geq -1} \lambda_q^{2s}\int_{\T^3}\Delta_q ((\nabla \times b)\times m)\cdot (\nabla \times m_q) \mathrm{d}x\\
=:& I + J.
\end{split}
\end{equation}

One further decomposes the terms $I$ and $J$ using Bony's paraproduct.
\begin{equation}\notag
\begin{split}
I = & \eta \sum_{q \geq -1}\sum_{|p-q| \leq 2} \lambda_q^{2s}\int_{\T^3}\Delta_q ((\nabla \times m_p)\times h_{\leq p-2})\cdot (\nabla \times m_q)\mathrm{d}x\\
& +\eta\sum_{q \geq -1}\sum_{|p-q| \leq 2} \lambda_q^{2s}\int_{\T^3} \Delta_q ((\nabla \times m_{\leq p-2})\times h_p)\cdot (\nabla \times m_q)\mathrm{d}x\\
& +\eta \sum_{q \geq -1}\sum_{p \geq q-2} \lambda_q^{2s}\int_{\T^3} \Delta_q ((\nabla \times \tilde m_{p})\times h_p)\cdot (\nabla \times m_q) \mathrm{d}x\\
=: & I_1+I_2+I_3;
\end{split}
\end{equation}
\begin{equation}\notag
\begin{split}
J = &\eta \sum_{q \geq -1} \sum_{|p-q| \leq 2} \lambda_q^{2s}\int_{\T^3} \Delta_q (m_{\leq p-2}\times (\nabla \times b_p))\cdot (\nabla \times m_q) \mathrm{d}x\\
& +\eta \sum_{q \geq -1} \sum_{|p-q| \leq 2}  \lambda_q^{2s}\int_{\T^3}\Delta_q (m_p\times(\nabla \times b_{\leq p-2}) )\cdot (\nabla \times m_q) \mathrm{d}x\\
& +\eta \sum_{q \geq -1} \sum_{p \geq q-2} \lambda_q^{2s}\int_{\T^3} \Delta_q (\tilde m_p\times(\nabla \times b_p) )\cdot (\nabla \times m_q) \mathrm{d}x\\
=:& J_1+J_2+J_3.
\end{split}
\end{equation}

One then proceeds to estimate the terms $I_1, I_2, I_3$ and $J_1, J_2, J_3.$ As for $I_1,$ one rewrites it using the commutator notation (\ref{cm1}) and notices that $I_{12}$ in the following expression vanishes.
\begin{equation}\notag
\begin{split}
I_1 = & \eta \sum_{q \geq -1}\sum_{|p-q| \leq 2}\lambda_q^{2s}\int_{\T^3}\big([\Delta_q, h_{\leq p-2} \times \nabla \times]m_p\big) \cdot (\nabla \times m_q) \mathrm{d}x\\
&-\eta \sum_{q \geq -1}\lambda_q^{2s}\int_{\T^3}\big(h_{\leq q-2}\times (\nabla \times m_q)\big)\cdot (\nabla \times m_q) \mathrm{d}x\\
&+ \eta \sum_{q \geq -1}\sum_{|p-q| \leq 2}\lambda_q^{2s}\int_{\T^3} \big( (h_{\leq q-2}-h_{\leq p-2})\times (\nabla \times (m_p)_q)\big)\cdot (\nabla \times m_q) \mathrm{d}x\\
=:& I_{11}+I_{12}+I_{13}.
\end{split}
\end{equation}

Taking into account that $m_{\leq Q_{b,h}}=0,$ one splits $I_{11}$ by the wavenumber.  
\begin{equation}\notag
\begin{split}
I_{11} = & \eta \sum_{q> Q_{b,h}}\sum_{|p-q| \leq 2}\lambda_q^{2s}\int_{\T^3}\big([\Delta_q, h_{\leq Q_{b,h}} \times \nabla \times]m_p\big) \cdot (\nabla \times m_q)\mathrm{d}x\\
&+\eta\sum_{q> Q_{b,h}}\sum_{|p-q| \leq 2}\lambda_q^{2s}\int_{\T^3}\big([\Delta_q, h_{(Q_h, p-2]}\times \nabla \times] m_p\big)\cdot (\nabla \times m_q) \mathrm{d}x\\
=:& I_{111}+I_{112}.
\end{split}
\end{equation}

By Lemma \ref{cmest}, H\"older's inequality, Definition \ref{wvnbr}, Young's inequalities, one estimates $I_{111}$ as follows.
\begin{equation}\notag
\begin{split}
|I_{111}| \leq & \eta \|\nabla h_{\leq Q_{b,h}}\|_\infty \sum_{q> Q_h}\lambda_q^{2s}\|m_q\|_2\sum_{|p-q| \leq 2}\|\nabla \times m_p\|_2\\
\leq & \eta \| h_{\leq Q_{b,h}}\|_\infty \sum_{q> Q_h}\lambda_q^{2s+1}\|m_q\|_2\sum_{|p-q| \leq 2}\|\nabla \times m_p\|_2\\
\lesssim & c_{r}\mu \sum_{q \geq -1}\lambda_q^{2s} \|\nabla m_q\|_2^2.
\end{split}
\end{equation}
One estimates $I_{112}$ using Lemma \ref{cmest}, H\"older's inequality, Definition \ref{wvnbr}, Young's and Jensen's inequalities.
\begin{equation}\notag
\begin{split}
|I_{112}| \leq & \eta \sum_{q> Q_h}\lambda_q^{2s}\|\nabla \times m_q\|_2\sum_{|p-q| \leq 2}\|m_p\|_2 \sum_{Q_h< p' \leq p-2}\lambda_{p'}\|h_{p'}\|_\infty \\
\leq & \eta \sum_{q> Q_h}\lambda_q^{2s}\|\nabla \times m_q\|_2\sum_{|p-q| \leq 2}\|\nabla \times m_p\|_2 \sum_{Q_h< p' \leq p-2}\| h_{p'}\|_\infty \lambda_{p'-p}\\
\leq & c_r {\mu} \sum_{q> Q_h}\lambda_q^{2s}\|\nabla \times m_q\|_2\sum_{|p-q| \leq 2}\|\nabla \times m_p\|_2 \sum_{Q_h< p' \leq q}\lambda_{q-p'}^{-1}\\
\lesssim & c_r {\mu} \sum_{q \geq -1}\lambda_q^{2s}\|\nabla m_q\|_2^2.
\end{split}
\end{equation}

For $p,q \in \mathbb{Z}$ satisfying $|p-q| \leq 2,$ it is true that $|h_{\leq q-2} - h _{\leq p-2}| \leq \sum_{i =0}^3|h_{q-i}|.$ Since $m_{q} = 0, \forall q \leq Q_{b,h},$ the following generic bound is true - 
$$ |I_{13}| \lesssim \eta \sum_{q > Q_{b,h}}\sum_{|p-q| \leq 2}\lambda_q^{2s}\int_{\T^3} |h_{q-2}| |\nabla \times (m_p)_q||\nabla \times m_q|\mathrm{d}x.\\ $$
The sum is then split by the wavenumber $Q_{b,h}.$
\begin{equation}\notag
\begin{split}
|I_{13}| \lesssim & \eta \sum_{Q_{b,h} < q \leq Q_{b,h}+2}\sum_{|p-q| \leq 2}\lambda_q^{2s}\int_{\T^3} |h_{q-2}||\nabla \times (m_p)_q||\nabla \times m_q|\mathrm{d}x\\
&+ \eta \sum_{q > Q_{b,h}+2}\sum_{|p-q| \leq 2}\lambda_q^{2s}\int_{\T^3} |h_{q-2}||\nabla \times (m_p)_q||\nabla \times m_q|\mathrm{d}x\\
=:& I_{131}+I_{132}.
\end{split}
\end{equation}

$I_{131}$ is estimated as follows.
\begin{equation}\notag
\begin{split}
I_{131}\leq & \eta \|h_{\leq Q_{b,h}}\|_\infty \sum_{Q_{b,h} <  q \leq Q_{b,h}+2}\lambda_q^{2s}\|\nabla \times m_q\|_2\sum_{|p-q| \leq 2}\|\nabla \times m_p\|_2\\
\leq & c_r {\mu} \sum_{Q_{b,h} < q \leq Q_{b,h}+2}\lambda_q^{2s}\|\nabla \times m_q\|_2\sum_{|p-q| \leq 2}\|\nabla \times m_p\|_2\\
\lesssim & c_r {\mu} \sum_{Q_{b,h}-2 < q \leq Q_{b,h}+2}\lambda_q^{2s} \|\nabla m_q\|_2^2.
\end{split}
\end{equation}
$I_{132}$ is estimated with H\"older's inequality, Definition \ref{wvnbr} and Young's inequality.
\begin{equation}\notag
\begin{split}
I_{132}\leq & \eta \sum_{q> Q_{b,h}+2}\lambda_q^{2s}\|h_{q-2}\|_\infty \|\nabla \times m_q\|_2\sum_{|p-q| \leq 2}\|\nabla \times m_p\|_2\\
\leq & c_r {\mu} \sum_{q> Q_{b,h}+2}\lambda_q^{2s}\|\nabla \times m_q\|_2\sum_{|p-q| \leq 2}\|\nabla \times m_p\|_2\\
\lesssim & c_r {\mu}\sum_{q \geq -1}\lambda_q^{2s} \|\nabla m_q\|_2^2.
\end{split}
\end{equation}

As $m_{\leq Q_{b,h}}=0,$ it is perceivable that $I_2$ consists of only high frequency parts and can be written as follows.
\begin{equation}\notag
\begin{split}
I_2 = \eta \sum_{p > Q_{b,h}+2}\sum_{|p-q| \leq 2}\lambda_q^{2s}\int_{\T^3}\Delta_q \big( h_p\times (\nabla \times m_{(Q_{b,h},  p-2]})\big)\cdot (\nabla \times m_q) \mathrm{d}x.
\end{split}
\end{equation}
Let $\delta > s.$  H\"older's, inequality, Definition \ref{wvnbr}, Young's and Jensen's inequalities lead to 
\begin{equation}\notag
\begin{split}
|I_{2}| \leq & \eta \sum_{p > Q_{b,h}} \|h_p\|_\infty \sum_{|p-q| \leq 2}\lambda_q^{2s} \|\nabla \times m_q\|_2 \sum_{Q_h < p' \leq p-2}  \|\nabla \times m_{p'}\|_2\\
\leq & \eta \sum_{q > Q_{b,h}-2} \lambda_q^{2s} \|\nabla \times m_q\|_2 \sum_{|p-q| \leq 2}\sum_{Q_{b,h} < p' \leq p-2}  \lambda_{p'}^\delta\|\nabla \times m_{p'}\|_2\lambda_{p-Q_{b,h}}^\delta\|h_p\|_\infty \\
\leq & c_r {\mu}\sum_{q > Q_{b,h}-2} \lambda_q^{2s} \|\nabla \times m_q\|_2 \sum_{|p-q| \leq 2}\sum_{Q_{b,h} < p' \leq p-2}  \lambda_{p'}^\delta\|\nabla \times m_{p'}\|_2 \\
\leq & c_r {\mu} \sum_{q >Q_{b,h}-2} \lambda_q^{s} \|\nabla \times m_q\|_2\sum_{Q_{b,h} < p' \leq q}  \lambda_{p'}^{s}\|\nabla \times m_{p'}\|_2\lambda_{p'-q}^{\delta-s}\\
\lesssim & c_r {\mu} \sum_{q \geq -1}\lambda_q^{2s}\|\nabla m_q\|_2^2.
\end{split}
\end{equation}

$I_3$ is split into three terms as follows.
\begin{equation}\notag
\begin{split}
I_{3}=& \eta \sum_{Q_{b,h} < q \leq Q_{b,h}+2}\sum_{q-2 \leq p \leq Q_{b,h}} \lambda_q^{2s}\int_{\T^3}\Delta_q ( h_p \times (\nabla \times \tilde m_{p}))\cdot \nabla \times m_q \mathrm{d}x\\
&+ \eta \sum_{Q_{b,h} < q \leq Q_{b,h}+2}\sum_{p > Q_{b,h}} \lambda_q^{2s}\int_{\T^3}\Delta_q ( h_p \times (\nabla \times \tilde m_{p}))\cdot \nabla \times m_q \mathrm{d}x\\
&+ \eta \sum_{q > Q_{b,h}+2}\sum_{p \geq q-2} \lambda_q^{2s}\int_{\T^3}\Delta_q ( h_p \times (\nabla \times \tilde m_{p}))\cdot \nabla \times m_q \mathrm{d}x\\
=:& I_{31}+I_{32}+I_{33}.
\end{split}
\end{equation}

Invoking Definition \ref{wvnbr} and applying H\"older's, Young's and Jensen's inequalities, one can estimate $I_{31}, I_{32}$ and $I_{33}$ as follows. 
\begin{equation}\notag
\begin{split}
|I_{31}| \leq & \eta \sum_{Q_{b,h} < q \leq Q_{b,h}+2}\lambda_q^{2s}\|\nabla \times m_q\|_2 \sum_{q-2 \leq p \leq Q_{b,h}}\|h_p\|_\infty \|\nabla \times \tilde m_p\|_2 \\ 
\lesssim & c_r {\mu} \sum_{Q_{b,h} < q \leq Q_{b,h}+2} \lambda_q^{s}\|\nabla \times m_q\|_2 \sum_{q-2 \leq p \leq Q_{b,h}}\lambda_p^{s}\|\nabla \times m_p\|_2\lambda_{q-p}^s\\
\lesssim & c_r {\mu} \sum_{Q_{b,h}-1 \leq q \leq Q_{b,h}+2}\lambda_q^{2s}\|\nabla m_q\|_2^2,
\end{split}
\end{equation}
\begin{equation}\notag
\begin{split}
|I_{32}| \leq & \sum_{Q_{b,h} < q \leq Q_{b,h}+2} \lambda_q^{2s}\|\nabla \times m_q\|_2 \sum_{p>Q_{b,h}}\|h_p\|_\infty \|\nabla \times \tilde m_p\|_2 \\ 
\lesssim & c_r {\mu}\sum_{Q_{b,h} < q \leq Q_{b,h}+2} \lambda_q^{s}\|\nabla \times m_q\|_2 \sum_{p > Q_{b,h}}\lambda_p^{s}\|\nabla \times m_p\|_2\lambda_{q-p}^{s} \lambda_{Q_{b,h}-p}^\delta \\
\lesssim & c_r {\mu}\sum_{q \geq -1}\lambda_q^{2s}\|\nabla m_q\|_2^2,
\end{split}
\end{equation}
\begin{equation}\notag
\begin{split}
|I_{33}| \leq & \sum_{q > Q_{b,h}+2} \lambda_q^{2s}\|\nabla \times m_q\|_2 \sum_{p \geq q-2}\|h_p\|_\infty \|\nabla \times \tilde m_p\|_2 \\ 
\lesssim & c_r {\mu} \sum_{q > Q_{b,h}+2} \lambda_q^{s}\|\nabla \times m_q\|_2 \sum_{p \geq q-2}\lambda_p^{s}\|\nabla \times m_p\|_2\lambda_{q-p}^{s}\lambda_{Q_{b,h}-p}^\delta \\
\lesssim & c_r {\mu} \sum_{q \geq -1}\lambda_q^{2s}\|\nabla m_q\|_2^2.
\end{split}
\end{equation}
Thus, the estimation for $I$ is completed.

$J_1, J_2$ and $J_3$ remain to be estimated. One can write $J_1,$ whose low frequency parts vanish due to $m_{\leq Q_{b,h}}=0$, as 
\begin{equation}\notag
\begin{split}
J_1 = & \sum_{p > Q_{b,h}+2}\sum_{|p-q| \leq 2}\lambda_q^{2s}\int_{\T^3}\Delta_q \big(m_{(Q_{b,h}, p-2]}\times (\nabla \times b_p)\big)\cdot (\nabla \times m_q) \mathrm{d}x.
\end{split}
\end{equation}
Recalling Definition \ref{wvnbr}, one can estimate $J_1$ using H\"older's, Young's and Jensen's inequalities, provided that $\delta>s+1.$ 
\begin{equation}\notag
\begin{split}
|J_1| \leq  & \sum_{p > Q_{b,h}+2}\lambda_p \|b_p\|_\infty \sum_{|p-q| \leq 2}\lambda_q^{2s} \| \nabla \times m_q\|_2\sum_{Q_{b,h} < p' \leq p-2 }\|m_{p'}\|_2\\
\leq  & \sum_{q > Q_{b,h}}\lambda_q^{2s} \| \nabla \times m_q\|_2\sum_{|p-q| \leq 2}\sum_{Q_{b,h} < p' \leq p-2 }\lambda_{p'}^\delta\|m_{p'}\|_2\lambda_{p-Q_{b,h}}^\delta  \|b_p\|_\infty\lambda_p^{1-\delta} \\
\leq  & c_r {\mu}\sum_{q > Q_{b,h}}\lambda_q^{s} \| \nabla \times m_q\|_2\sum_{Q_{b,h} < p' \leq q }\lambda_{p'}^{s+1}\|m_{p'}\|_2\lambda_{p'-q}^{\delta-s-1}\\
\lesssim  & c_r {\mu} \sum_{q \geq -1}\lambda_q^{2s} \|\nabla m_q\|_2^2.
\end{split}
\end{equation}

$J_2$ can be partitioned into two terms with $Q_{b,h}$.
\begin{equation}\notag
\begin{split}
J_2 = & \sum_{q > Q_{b,h}} \sum_{|p-q| \leq 2} \lambda_q^{2s}\int_{\T^3}\Delta_q \big(m_p\times(\nabla \times b_{\leq Q_{b,h}}) \big)\cdot (\nabla \times m_q) \mathrm{d}x\\
& + \sum_{q > Q_{b,h}} \sum_{|p-q| \leq 2} \lambda_q^{2s}\int_{\T^3}\Delta_q \big(m_p\times(\nabla \times b_{(Q_{b,h}, p-2]})\big)\cdot (\nabla \times m_q) \mathrm{d}x\\
=:& J_{21}+J_{22}.
\end{split}
\end{equation}
To estimate $J_{21},$ one applies H\"older's and Young's inequalities.
\begin{equation}\notag
\begin{split}
|J_{21}| \leq & \eta \|\nabla b_{\leq Q_{b,h}}\|_\infty \sum_{q > Q_{b,h}}\lambda_q^{2s}\|\nabla \times m_q\|_2 \sum_{|p-q| \leq 2}\|m_p\|_2\\
\lesssim & \eta \|b_{\leq Q_{b,h}}\|_\infty \sum_{q > Q_{b,h}}\lambda_q^{2s+1}\|\nabla m_q\|_2 \sum_{|p-q| \leq 2}\|m_p\|_2\\
\lesssim & c_r {\mu} \sum_{q > Q_{b,h}}\lambda_q^{s}\|\nabla m_q\|_2 \sum_{|p-q| \leq 2}\lambda_p^{s+1}\|m_p\|_2\\
\lesssim  & c_r {\mu}\sum_{q \geq -1}\lambda_q^{2s} \|\nabla m_q\|_2^2.
\end{split}
\end{equation}
For $J_{22},$ H\"older's inequality, Definition \ref{wvnbr}, Young's and Jensen's inequalities yield 
\begin{equation}\notag
\begin{split}
|J_{22}| \leq &  \eta \sum_{q > Q_{b,h}}\lambda_q^{2s}\|\nabla \times m_q\|_2 \sum_{|p-q| \leq 2}\|m_p\|_2\sum_{Q_{b,h}< p' \leq p-2}\lambda_{p'}\|b_{p'}\|_\infty\\
\leq &  c_r {\mu} \sum_{q > Q_{b,h}}\lambda_q^{s}\|\nabla \times m_q\|_2 \sum_{|p-q| \leq 2}\lambda_p^{s+1}\|m_p\|_2\sum_{Q_{b,h}< p' \leq p-2}\lambda_{p'-p}\\
\lesssim  & c_r {\mu} \sum_{q \geq -1}\lambda_q^{2s} \|\nabla m_q\|_2^2.
\end{split}
\end{equation}

Taking advantage of $m_{\leq Q_{b,h}}=0,$ one write $J_3$ as  
\begin{equation}\notag
\begin{split}
|J_3|=& \eta \sum_{q > Q_{b,h}} \sum_{p \geq q+2} \lambda_q^{2s}\int_{\T^3}\Delta_q \big(\tilde m_p\times(\nabla \times b_p) \big)\cdot (\nabla \times m_q)\mathrm{d}x,
\end{split}
\end{equation}
which can then be estimated as follows.
\begin{equation}\notag
\begin{split}
|J_3| \leq & \eta\sum_{p \geq Q_{b,h}} \lambda_p\|b_p\|_\infty \|m_p\|_2  \sum_{q \leq p-2}\lambda_q^{2s} \|\nabla \times m_q\|_2 \\
\leq & c_r {\mu} \sum_{q > Q_{b,h}-2}\lambda_q^{2s} \|\nabla \times m_q\|_2 \sum_{p \geq q+2} \lambda_p\|m_p\|_2\\
\leq & c_r {\mu} \sum_{q > Q_{b,h}-2}\lambda_q^{s} \|\nabla \times m_q\|_2 \sum_{p \geq q+2} \lambda_p^{s+1}\|m_p\|_2\lambda_{q-p}^{s}\lambda_{Q_{b,h}-p}^\delta\\
\lesssim & c_r {\mu} \sum_{q \geq -1}\lambda_q^{2s} \|\nabla m_q\|_2^2.
\end{split}
\end{equation}

Let $c_r = 1-(2\mu)^{-1}.$ Assembling all the estimates above leads to 
\begin{equation}\notag
\frac{\mathrm{d}}{\mathrm{d}t}\sum_{q \geq -1}\lambda_q^{2s}\|m_q\|_2^2 \lesssim -\sum_{q \geq -1}\lambda_q^{2s+2}\|m_q\|_2^2 \lesssim - \lambda_0^2 \sum_{q \geq -1}\lambda_q^{2s}\|m_q\|_2^2.
\end{equation}
Setting $s=0,$ one sees that
$
\frac{\mathrm{d}}{\mathrm{d}t}\|m\|_{L^2}^2 \lesssim -\lambda_0^2 \|m\|_{L^2}^2.
$
The desired result  then follows from Gr\"onwall's inequality.

\cbdu

\section{Proof of Theorem \ref{thm}}

Let $(u,b)$ and $(v,h)$ be two weak solutions to system (\ref{hmhd1})-(\ref{hmhd3}). Straightforward calculations show that the difference $(w, m):=(u-v, b-h)$ satisfies the following system of equations.
\begin{equation}\label{dfrnc}
\begin{split}
w_t-\nu\Delta w=& -u\cdot \nabla w-w \cdot \nabla v+b \cdot \nabla m+m\cdot \nabla h-\nabla \pi,\\
m_t-\mu\Delta m=& -v\cdot \nabla m-w\cdot \nabla b+b\cdot \nabla w+m\cdot \nabla v\\
&-\nabla \times (\nabla \times m)\times h)-\nabla \times ((\nabla \times b)\times m).\\
\end{split}
\end{equation}

Utilizing the wavenumbers, one shall eventually prove that $(w,m)$ satisfies the following inequality
\begin{equation}
\frac{\mathrm{d}}{\mathrm{d}t}\big(\|w\|_{L^2}^2+\|m\|_{L^2}^2\big) \lesssim -\big(\|\nabla w \|_{L^2}^2 + \|\nabla m\|_{L^2}^2\big),
\end{equation}
which leads to theorem (\ref{thm}).

To this end, one considers a frequency-localized version of system (\ref{dfrnc}) in energy spaces. Multiplying the equations by $\lambda_q^{2s}\Delta_q^2 w$ and $\lambda_q^{2s}\Delta_q^2 m$ respectively, integrating by parts and summing over $q,$ one obtains
\begin{equation}
\begin{split}
&\frac{1}{2}\frac{\mathrm{d}}{\mathrm{d}t}\sum_{q \geq -1}\lambda_q^{2s}\|w_q\|_2^2+\nu\sum_{q \geq -1}\lambda_q^{2s+2}\|w_q\|_2^2\\
\leq & -\sum_{q \geq -1} \lambda_q^{2s}\int_{\T^3} \Delta_q(u\cdot \nabla w)\cdot w_q \mathrm{d}x -\sum_{q \geq -1} \lambda_q^{2s}\int_{\T^3} \Delta_q(w\cdot \nabla v)\cdot w_q \mathrm{d}x\\
&+\sum_{q \geq -1} \lambda_q^{2s}\int_{\T^3}\Delta_q(b\cdot \nabla m)\cdot w_q \mathrm{d}x +\sum_{q \geq -1} \lambda_q^{2s}\int_{\T^3} \Delta_q(m\cdot \nabla h)\cdot w_q \mathrm{d}x\\
=:&A+B+C+D,\\
\end{split}
\end{equation}
and
\begin{equation}
\begin{split}
&\frac{1}{2}\frac{\mathrm{d}}{\mathrm{d}t}\sum_{q \geq -1}\lambda_q^{2s}\|m_q\|_2^2+\mu\sum_{q \geq -1}\lambda_q^{2s+2}\|m_q\|_2^2\\
\leq & -\sum_{q \geq -1} \lambda_q^{2s}\int_{\T^3}  \Delta_q(v\cdot \nabla m)\cdot m_q \mathrm{d}x -\lambda_q^{2s}\int_{\T^3} \Delta_q(w\cdot \nabla b) \cdot m_q \mathrm{d}x\\
&+\sum_{q \geq -1} \lambda_q^{2s}\int_{\T^3} \Delta_q(b\cdot \nabla w)\cdot m_q \mathrm{d}x +\sum_{q \geq -1} \lambda_q^{2s}\int_{\T^3} \Delta_q(m\cdot \nabla v)\cdot m_q \mathrm{d}x\\
&- \sum_{q \geq -1} \lambda_q^{2s}\int_{\T^3}\Delta_q ( (\nabla \times m)\times h)\cdot (\nabla \times m_q)\mathrm{d}x\\
&- \sum_{q \geq -1} \lambda_q^{2s}\int_{\T^3}\Delta_q ((\nabla \times b)\times m)\cdot (\nabla \times m_q) \mathrm{d}x\\
=:& E+F+G+H+I+J.
\end{split}
\end{equation}
The tasks are then to control the terms $A$--$J.$

\subsection{Estimation of A}

The estimates for $A$ fall into the same line as those in \cite{CDK}. Bony's decomposition leads to the following - 
\begin{equation}\notag
\begin{split}
A=&-\sum_{q \geq -1}\sum_{|p-q|\leq 2} \lambda_q^{2s}\int_{\T^3} \Delta_q(u_{\leq p-2}\cdot \nabla w_p)\cdot w_q \mathrm{d}x\\
&-\sum_{q \geq -1}\sum_{|p-q|\leq 2} \lambda_q^{2s}\int_{\T^3} \Delta_q(u_{p}\cdot \nabla w_{\leq p-2})\cdot w_q \mathrm{d}x\\
&-\sum_{q \geq -1}\sum_{p \geq q-2} \lambda_q^{2s}\int_{\T^3} \Delta_q(u_p\cdot \nabla \tilde w_p)\cdot w_q \mathrm{d}x\\
=: & A_1+A_2+A_3.
\end{split}
\end{equation}

Using Definition \ref{uvnbr}, one then separate the lower and higher modes of $A_1.$
\begin{equation}\notag
\begin{split}
|A_1| \leq & \sum_{p > Q_{u,v}}\sum_{|q-p|\leq 2} \lambda_q^{2s}\int_{\T^3} |\Delta_q(u_{\leq p-2}\cdot \nabla w_p)\cdot w_q| \mathrm{d}x\\
\leq & \sum_{p > Q_{u,v}}\sum_{|q-p|\leq 2} \lambda_q^{2s}\int_{\T^3} |\Delta_q(u_{\leq Q_{u,v}}\cdot \nabla w_p)\cdot w_q| \mathrm{d}x\\
& + \sum_{p' > Q_{u,v}} \sum _{p \geq p'+2}\sum_{|q-p|\leq 2} \lambda_q^{2s}\int_{\T^3} |\Delta_q(u_{p'}\cdot \nabla w_p)\cdot w_q| \mathrm{d}x\\
=: & A_{11} + A_{12}.
\end{split}
\end{equation}

To control the lower modes, one uses Definition \ref{uvnbr}, Lemma \ref{brn}, H\"older's and Young's inequalities.
\begin{equation}\notag
\begin{split}
A_{11} \lesssim &   \|u_{\leq Q_{u,v}}\|_r  \sum_{p > Q_{u,v}}\lambda_p\|w_p\|_2\sum_{|q-p|\leq 2} \lambda_q^{2s} \|w_q\|_\frac{2r}{r-2}\\
\lesssim &  c_{r}\nu \sum_{p > Q_{u,v}}\lambda_p\|w_p\|_2\sum_{|q-p|\leq 2} \lambda_{Q_{u,v}}^{1-\frac{3}{r}}\lambda_q^{2s+\frac{3}{r}}\|w_q\|_2\\
\lesssim &  c_{r}\nu \sum_{q \geq -1}\lambda_q^{2s} \| \nabla w_q\|_2^2.
\end{split}
\end{equation}

The higher modes are estimated as follows.
\begin{equation}\notag
\begin{split}
A_{12} \lesssim &   \sum_{p' > Q_{u,v}} \|u_{p'}\|_r  \sum_{p > p'+2}\lambda_p\|w_p\|_2\sum_{|q-p|\leq 2} \lambda_q^{2s}\|w_q\|_\frac{2r}{r-2}\\
\lesssim &  c_{r}\nu \sum_{p' > Q_{u,v}} \lambda_{p'}^{1-\frac{3}{r}} \sum_{p > p'+2}\lambda_p\|w_p\|_2\sum_{|q-p|\leq 2} \lambda_q^{2s+\frac{3}{r}}\|w_q\|_2\\
\lesssim &  c_{r}\nu \sum_{q \geq -1}\lambda_q^{2s} \| \nabla w_q\|_2^2.
\end{split}
\end{equation}

It follows from the condition $w_{\leq Q_{u,v}}=0$ that
\begin{equation}\notag
\begin{split}
A_2 = &-\sum_{p > Q_{u,v}+2}\sum_{|q-p|\leq 2} \lambda_q^{2s}\int_{\T^3} \Delta_q(u_{p}\cdot \nabla w_{\leq p-2})\cdot w_q \mathrm{d}x.
\end{split}
\end{equation}

Recalling Definition \ref{uvnbr}, one then estimates $A_2$ using H\"older's and Young's inequalities.
\begin{equation}\notag
\begin{split}
|A_2| \leq & \sum_{p > Q_{u,v}+2}\|u_p\|_r\sum_{Q_{u,v} < p' \leq p+2}\lambda_{p'}\|w_{p'}\|_{\frac{2r}{r-2}}\sum_{|q-p|\leq 2} \lambda_q^{2s}\|w_q\|_2\\
\lesssim &  c_{r}\nu\sum_{p > Q_{u,v}+2}\lambda_p^{1-\frac{3}{r}}\sum_{Q_{u,v} < p' \leq p+2}\lambda_{p'}^{1+\frac{3}{r}}\|w_{p'}\|_2\sum_{|q-p|\leq 2} \lambda_q^{2s}\|w_q\|_2\\
\lesssim &  c_{r}\nu\sum_{q > Q_{u,v}}\lambda_q^{s+1}\|w_q\|_2\sum_{Q_{u,v} < p' \leq q}\lambda_{p'}^{s+1}\|w_{p'}\|_2\lambda_{p'-q}^{\frac{3}{r}-s}\\
\lesssim &  c_{r}\nu \sum_{q \geq -1}\lambda_q^{2s} \| \nabla w_q\|_2^2.
\end{split}
\end{equation}

Separating lower and higher modes of $A_3$ with the wavenumber $Q_{u,v}$ results in
\begin{equation}\notag
\begin{split}
A_3 = &-\sum_{p= Q_{u,v}}\sum_{q \leq p+2 } \lambda_q^{2s}\int_{\T^3} \Delta_q(u_p\cdot \nabla \tilde w_p)\cdot w_q \mathrm{d}x\\
 &-\sum_{p > Q_{u,v}}\sum_{q \leq p+2 } \lambda_q^{2s}\int_{\T^3} \Delta_q(u_p\cdot \nabla \tilde w_p)\cdot w_q \mathrm{d}x\\
 =: & A_{31}+A_{32}.
\end{split}
\end{equation}

One has no difficulty in controlling the few lower modes.
\begin{equation}\notag
\begin{split}
|A_{31}| \lesssim & \Lambda_{u,v}\|u_{Q_{u,v}}\|_r\|w_{Q_{u,v}}\|_2\sum_{Q_{u,v} < q \leq Q_{u,v}+2 }\lambda_q^{2s}\|w_q\|_{\frac{2r}{r-2}}\\
\lesssim & c_{r} \nu  \Lambda_{u,v}^{2-\frac{3}{r}}\|w_{Q_{u,v}}\|_2\sum_{Q_{u,v} < q \leq Q_{u,v}+2 }\lambda_q^{2s+\frac{3}{r}}\|w_q\|_2\\
\lesssim & c_{r} \nu  \sum_{Q_{u,v} < q \leq Q_{u,v}+2 }\lambda_q^{2s}\|\nabla w_{q}\|_2^2.
\end{split}
\end{equation}

The higher modes are estimated using Definition \ref{uvnbr}, H\"older's, Young's and Jensen's inequalities.
\begin{equation}\notag
\begin{split}
|A_{32}| \leq & \sum_{p > Q_{u,v}}\|u_p\|_r\|\nabla \tilde w_p\|_2\sum_{q \leq p+2 } \lambda_q^{2s}\|w_q\|_{\frac{2r}{r-2}}\\
\lesssim & c_{r} \nu \sum_{p > Q_{u,v}}\lambda_p^{2-\frac{3}{r}}\|w_p\|_2\sum_{q \leq p+2 } \lambda_q^{2s+\frac{3}{r}}\|w_q\|_{2}\\
\lesssim & c_{r} \nu \sum_{p > Q_{u,v}}\lambda_p^{s+1}\|w_p\|_2\sum_{q \leq p+2 } \lambda_q^{s+1}\|w_q\|_2\lambda_{q-p}^{s-1+\frac{3}{r}}\\
\lesssim &  c_{r} \nu \sum_{q \geq -1}\lambda_q^{2s} \| \nabla w_q\|_2^2.
\end{split}
\end{equation}

\subsection{Estimation of B}

As a result of Bony's paraproduct decomposition
\begin{equation}\notag
\begin{split}
B =& -\sum_{q \geq -1} \sum_{|p-q|\leq 2} \lambda_q^{2s}\int_{\T^3} \Delta_q(w_{\leq p-2}\cdot \nabla v_p)\cdot w_q \mathrm{d}x\\
&-\sum_{q \geq -1} \sum_{|p-q|\leq 2} \lambda_q^{2s}\int_{\T^3} \Delta_q(w_{p}\cdot \nabla v_{ \leq p-2})\cdot w_q \mathrm{d}x\\
&-\sum_{q \geq -1} \sum_{p \geq q-2} \lambda_q^{2s}\int_{\T^3} \Delta_q(w_{p}\cdot \nabla \tilde v_{p})\cdot w_q \mathrm{d}x\\
=:& B_1+B_2+B_3.
\end{split}
\end{equation}

Since $w_{\leq Q_{u,v}} =0,$ $B_1$ consists of only higher modes.
\begin{equation}\notag
B_1 = -\sum_{p> Q_{u,v}+2} \sum_{|q-p|\leq 2} \lambda_q^{2s}\int_{\T^3} \Delta_q(w_{\leq p-2}\cdot \nabla v_p)\cdot w_q \mathrm{d}x.
\end{equation}

Let $1-\frac{3}{r} <0.$ One can estimate $B_1$ using Definition \ref{uvnbr}, H\"older's, Young's and Jensen's inequalities. 
\begin{equation}\notag
\begin{split}
|B_1| \lesssim & \sum_{p> Q_{u,v}+2}\lambda_p \|v_p\|_r \sum_{Q_{u,v} < p' \leq p-2}\|w_{p'}\|_{\frac{2r}{r-2}} \sum_{|q-p|\leq 2} \lambda_q^{2s}\|w_q\|_2\\
\lesssim & c_{r} \nu \sum_{p> Q_{u,v}+2}\lambda_p^{2-\frac{3}{r}} \sum_{Q_{u,v} < p' \leq p-2}\lambda_{p'}^\frac{3}{r}\|w_{p'}\|_2 \sum_{|q-p|\leq 2} \lambda_q^{2s}\|w_q\|_2\\
\lesssim & c_{r} \nu \sum_{p> Q_{u,v}+2}\lambda_p^{2s+2-\frac{3}{r}}\|w_p\|_2 \sum_{Q_{u,v} < p' \leq p-2}\lambda_{p'}^\frac{3}{r}\|w_{p'}\|_2\\
\lesssim & c_{r} \nu \sum_{p> Q_{u,v}+2}\lambda_p^{s+1}\|w_p\|_2 \sum_{Q_{u,v} < p' \leq p-2}\lambda_{p'}^{s+1}\|w_{p'}\|_2\lambda_{p-p'}^{s+1-\frac{3}{r}}\\
\lesssim &  c_{r} \nu  \sum_{q \geq -1}\lambda_q^{2s} \| \nabla w_q\|_2^2.
\end{split}
\end{equation}

Splitting $B_2$ with the wavenumber $Q_{u,v}$ results in
\begin{equation}\notag
\begin{split}
B_2 = &-\sum_{Q_{u,v}< p \leq Q_{u,v}+2} \sum_{|q-p|\leq 2} \lambda_q^{2s}\int_{\T^3} \Delta_q(w_{p}\cdot \nabla v_{ \leq p-2})\cdot w_q \mathrm{d}x\\
&-\sum_{p > Q_{u,v}+2} \sum_{|q-p|\leq 2} \lambda_q^{2s}\int_{\T^3} \Delta_q(w_{p}\cdot \nabla v_{ \leq Q_{u,v}})\cdot w_q \mathrm{d}x\\
&-\sum_{p > Q_{u,v}+2} \sum_{|q-p|\leq 2} \lambda_q^{2s}\int_{\T^3} \Delta_q(w_{p}\cdot \nabla v_{(Q_{u,v}, p-2]})\cdot w_q \mathrm{d}x\\
=:& B_{21}+B_{22}+B_{23}.
\end{split}
\end{equation}

The estimate for the lower modes $|B_{21}|+|B_{22}|$ are as follows.
\begin{equation}\notag
\begin{split}
|B_{21}|+|B_{22}| \lesssim & \sum_{p>Q_{u,v}}\|w_p\|_{\frac{2r}{r-2}}\sum_{|q-p|\leq 2}\lambda_q^{2s}\|w_q\|_{2}\sum_{p'<Q_{u,v}}\lambda_{p'}\|v_{p'}\|_r\\
\lesssim & c_{r} \nu \sum_{p>Q_{u,v}}\lambda_p^{\frac{3}{r}}\|w_p\|_2\sum_{|q-p|\leq 2}\lambda_q^{2s}\|w_q\|_{2}\sum_{p'<Q_{u,v}}\lambda_{p'}^{1-\frac{3}{r}}\\
\lesssim & c_{r} \nu \sum_{p>Q_{u,v}}\lambda_p^2\|w_p\|_2\sum_{|q-p|\leq 2}\lambda_q^{2s}\|w_q\|_{2}\sum_{p'<Q_{u,v}}\lambda_{p'-p}^{2-\frac{3}{r}}\\
\lesssim & c_{r} \nu  \sum_{q \geq -1}\lambda_q^{2s} \| \nabla w_q\|_2^2.
\end{split}
\end{equation}

The estimate for $B_{23}$ follows from Definition \ref{uvnbr} and H\"older's inequality.
\begin{equation}\notag
\begin{split}
|B_{23}| \leq & \sum_{p > Q_{u,v}+2}\|w_p\|_2\sum_{Q_{u,v}< p' \leq p-2}\|v_{p'}\|_r \sum_{|q-p|\leq 2} \lambda_q^{2s}\|w_q\|_{\frac{2r}{r-2}}\\
\leq & c_{r} \nu \sum_{p > Q_{u,v}+2}\|w_p\|_2\sum_{Q_{u,v}< p' \leq p-2}\lambda_{p'}^{2-\frac{3}{r}} \sum_{|q-p|\leq 2} \lambda_q^{2s+\frac{3}{r}}\|w_q\|_2\\
\leq & c_{r} \nu \sum_{p > Q_{u,v}+2}\lambda_p^{s+1}\|w_p\|_2\sum_{|q-p|\leq 2}\lambda_q^{s+1}\|w_q\|_2\sum_{Q_{u,v}< p' \leq p-2}\lambda_{p'-p}^{2-\frac{3}{r}} \\
\lesssim &  c_{r} \nu  \sum_{q \geq -1}\lambda_q^{2s} \| \nabla w_q\|_2^2.
\end{split}
\end{equation}

Similar to previous terms, $|B_3|$ is bounded above by the estimates for the lower modes and for the higher modes.
\begin{equation}\notag
\begin{split}
|B_3|\leq & \sum_{p= Q_{u,v}+1} \sum_{q \leq p+2} \lambda_q^{2s}\int_{\T^3} |\Delta_q(w_{p}\cdot \nabla v_{p-1})\cdot w_q| \mathrm{d}x\\ 
&+\sum_{p> Q_{u,v}+1} \sum_{q \leq p+2} \lambda_q^{2s}\int_{\T^3} |\Delta_q(w_{p}\cdot \nabla v_{p-1})\cdot w_q| \mathrm{d}x\\
=: & B_{31}+B_{32}.
\end{split}
\end{equation}

The term $B_{31},$ consisting of scarce lower modes, can be controlled with ease.
\begin{equation}\notag
\begin{split}
B_{31} \lesssim & \Lambda_{u,v}^{2s+1}\|w_{Q_{u,v}}\|_2\|v_{Q_{u,v}}\|_r \sum_{Q_{u,v}<q \leq Q_{u,v}+3}\|w_q\|_\frac{2r}{r-2}\\
\lesssim & c_{r} \nu  \Lambda_{u,v}^{2s+2-\frac{3}{r}}\|w_{Q_{u,v}}\|_2 \sum_{Q_{u,v}<q \leq Q_{u,v}+3}\lambda_q^{\frac{3}{r}}\|w_q\|_2\\
\lesssim & c_{r} \nu  \sum_{Q_{u,v} < q \leq Q_{u,v}+3 }\lambda_q^{2s}\|\nabla w_{q}\|_2^2.
\end{split}
\end{equation}

As a result of Definition \ref{uvnbr}, H\"older's, Young's and Jensen's inequalities, $B_{32}$ can be estimated as follows.
\begin{equation}\notag
\begin{split}
B_{32} \lesssim & \sum_{p> Q_{u,v}+1}\lambda_p\|v_p\|_r\|w_p\|_2\sum_{q \leq p+2} \lambda_q^{2s}\|w_q\|_{\frac{2r}{r-2}}\\
\lesssim & c_{r} \nu \sum_{p> Q_{u,v}+1}\lambda_p^{2-\frac{3}{r}}\|w_p\|_2\sum_{q \leq p+2} \lambda_q^{2s+\frac{3}{r}}\|w_q\|_2\\
\lesssim & c_{r} \nu \sum_{p> Q_{u,v}+1}\lambda_p^{s+1}\|w_p\|_2\sum_{q \leq p+2} \lambda_q^{s+1}\|w_q\|_2\lambda_{q-p}^{s-1+\frac{3}{r}}\\
\lesssim &  c_{r} \nu  \sum_{q \geq -1}\lambda_q^{2s} \| \nabla w_q\|_2^2.
\end{split}
\end{equation}

\subsection{Estimation of C}

Bony's paraproduct decomposition yields
\begin{equation}\notag
\begin{split}
C = & \sum_{q \geq -1} \sum_{|p-q| \leq 2} \lambda_q^{2s}\int_{\T^3}\Delta_q(b_{\leq p-2}\cdot \nabla m_p)w_q \mathrm{d}x\\
&+ \sum_{q \geq -1} \sum_{|p-q| \leq 2} \lambda_q^{2s}\int_{\T^3}\Delta_q(b_p \cdot \nabla m_{\leq p-2})w_q \mathrm{d}x\\
&+ \sum_{q \geq -1} \sum_{p \geq q-2} \lambda_q^{2s}\int_{\T^3}\Delta_q(b_p \cdot \nabla \tilde m_p)w_q \mathrm{d}x\\
=:& C_1+C_2+C_3.
\end{split}
\end{equation}

Moreover, one rewrites $C_1$ using the commutator as
\begin{equation}\notag
\begin{split}
C_1 =& \sum_{q \geq -1} \sum_{|p-q| \leq 2} \lambda_q^{2s}\int_{\T^3}[\Delta_q, b_{\leq p-2}\cdot \nabla] m_p w_q \mathrm{d}x\\
&+ \sum_{q \geq -1} \sum_{|p-q| \leq 2} \lambda_q^{2s}\int_{\T^3}b_{\leq q-2} \cdot \nabla \Delta_q m_p w_q \mathrm{d}x\\
&+ \sum_{q \geq -1}\sum_{|p-q| \leq 2} \lambda_q^{2s}\int_{\T^3}(b_{\leq p-2}-b_{\leq q-2})\cdot \nabla \Delta_q m_p w_q \mathrm{d}x\\
=:& C_{11}+C_{12}+C_{13}.
\end{split}
\end{equation}
As will be seen later, $C_{12}$ cancels a part of the term $G.$ 

Taking into account that $m_{\leq Q_{b,h}}=0,$ one split $C_{11}$ using the wavenumber $Q_{b,h}.$
\begin{equation}\notag
\begin{split}
C_{11} = & \sum_{Q_{b,h}< p \leq Q_{b,h}+2} \sum_{|q-p| \leq 2} \lambda_q^{2s}\int_{\T^3}[\Delta_q, b_{\leq p-2}\cdot \nabla] m_p w_q \mathrm{d}x\\
&+ \sum_{p > Q_{b,h}+2} \sum_{|q-p| \leq 2} \lambda_q^{2s}\int_{\T^3}[\Delta_q, b_{\leq Q_{b,h}}\cdot \nabla] m_p w_q \mathrm{d}x\\
& + \sum_{p > Q_{b,h}+2} \sum_{|q-p| \leq 2} \lambda_q^{2s}\int_{\T^3}[\Delta_q, b_{(Q_{b,h}, p-2]}\cdot \nabla] m_p w_q \mathrm{d}x\\
=:& C_{111}+C_{112}+C_{113}.
\end{split}
\end{equation}

By Definition \ref{wvnbr}, H\"older's and Young's inequalities, the following estimate holds.
\begin{equation}\notag
\begin{split}
|C_{111}|+|C_{112}| \leq & \|\nabla b_{\leq Q_{b,h}} \|_\infty \sum_{p > Q_{b,h}}\|m_p\|_2 \sum_{|q-p|\leq 2}\lambda_q^{2s}\| w_q\|_2\\
\leq & c_r \kappa \sum_{p \geq -1}\lambda_p^{s+1}\|m_p\|_2 \sum_{|q-p|\leq 2}\lambda_q^{s+1}\| w_q\|_2\\
\leq & c_r \kappa \sum_{q \geq -1 }(\lambda_q^{2s+2}\|w_q\|_2^2+\lambda_q^{2s+2}\|m_q\|_2^2).
\end{split}
\end{equation}

As a result of Definition \ref{wvnbr}, H\"older's and Young's inequalities, the following estimate for $C_{113}$ is true.
\begin{equation}\notag
\begin{split}
|C_{113}| \leq & \sum_{p > Q_{b,h}}\|m_p\|_2 \sum_{|q-p|\leq 2}\lambda_q^{2s}\| w_q\|_2 \sum_{Q_{b,h}< p' \leq p-2}\lambda_{p'}\|b_{p'}\|_\infty\\
\leq & c_r \kappa \sum_{p > Q_{b,h}}\|m_p\|_2 \sum_{|q-p|\leq 2}\lambda_q^{2s+2}\| w_q\|_2 \sum_{Q_{b,h}< p' \leq p-2}\lambda_{p'-p}^2\\
\leq & c_r \kappa \sum_{p > Q_{b,h}}\lambda_p^{s+1}\|m_p\|_2 \sum_{|q-p|\leq 2}\lambda_q^{s+1}\| w_q\|_2 \sum_{Q_{b,h}< p' \leq p-2}\lambda_{p'-p}^2\\
\leq & c_r \kappa \sum_{q \geq -1 }(\lambda_q^{2s+2}\|w_q\|_2^2+\lambda_q^{2s+2}\|m_q\|_2^2).
\end{split}
\end{equation}

$|C_{13}|$ is bounded above by two terms as follows.
\begin{equation}\notag
\begin{split}
|C_{13}| \lesssim & \sum_{q \geq -1}\sum_{|p-q| \leq 2} \lambda_q^{2s}\int_{\T^3}(|b_{q-3}|+|b_{q-2}|+|b_{q-1}|+|b_{q}| )|\nabla \Delta_q m_p w_q|\mathrm{d}x\\
\lesssim & \sum_{-1 \leq q \leq Q_{b,h}}\sum_{|p-q| \leq 2} \lambda_q^{2s}\int_{\T^3}|b_{q}||\nabla \Delta_q m_p w_q|\mathrm{d}x\\ &+\sum_{q > Q_{b,h}}\sum_{|p-q| \leq 2} \lambda_q^{2s}\int_{\T^3}|b_{q}||\nabla \Delta_q m_p w_q|\mathrm{d}x\\
=:& C_{131}+C_{132}.
\end{split}
\end{equation}

The estimate for $C_{131}$ is as follows. 
\begin{equation}\notag
\begin{split}
C_{131} \leq & \sum_{-1 \leq q \leq Q_{b,h}}\|b_q\|_\infty \lambda_q^{2s}\|w_q\|_2 \sum_{|p-q| \leq 2} \|\nabla m_p\|_2\\
\lesssim & c_r \kappa \sum_{q \geq -1}\lambda_q^{s+1}\|w_q\|_2 \sum_{|p-q| \leq 2} \lambda_p^{s+1}\|m_p\|_2\\
\lesssim & c_r \kappa \sum_{q \geq -1 }(\lambda_q^{2s+2}\|w_q\|_2^2+\lambda_q^{2s+2}\|m_q\|_2^2)
\end{split}
\end{equation}

$C_{132}$ enjoys the following estimate, thanks to Definition \ref{wvnbr}.
\begin{equation}\notag
\begin{split}
C_{132} \leq & \sum_{q > Q_{b,h}}\|b_q\|_\infty \lambda_q^{2s}\|w_q\|_2 \sum_{|p-q| \leq 2} \|\nabla m_p\|_2\\
\lesssim & c_r \kappa \sum_{q \geq -1}\lambda_q^{s+1}\|w_q\|_2 \sum_{|p-q| \leq 2} \lambda_p^{s+1}\|m_p\|_2\\
\lesssim & c_r \kappa \sum_{q \geq -1 }(\lambda_q^{2s+2}\|w_q\|_2^2+\lambda_q^{2s+2}\|m_q\|_2^2)
\end{split}
\end{equation}

Since $\nabla m_{\leq Q_{b,h}}=0,$ the lower modes of $C_2$ vanish and it can be seen that
\begin{equation}\notag
\begin{split}
C_2 =& \sum_{p > Q_{b,h}+2} \sum_{|p-q| \leq 2} \lambda_q^{2s}\int_{\T^3}\Delta_q(b_p \cdot \nabla m_{(Q_{b,h}, p-2]})w_q \mathrm{d}x,
\end{split}
\end{equation}
which is estimated using H\"older's, Young's and Jensen's inequalities as
\begin{equation}\notag
\begin{split}
|C_2| \leq & \sum_{p > Q_{b,h}+2}\|b_p\|_\infty \sum_{|p-q| \leq 2}\lambda_q^{2s}\|w_q\|_2 \sum_{Q_{b,h}< p' \leq p-2}\lambda_{p'} \|m_{p'}\|_2\\
\leq & \sum_{p > Q_{b,h}+2}\|b_p\|_\infty \sum_{|p-q| \leq 2}\lambda_q^{s+1}\|w_q\|_2 \sum_{Q_{b,h}< p' \leq p-2}\lambda_{p'}^{s+1} \|m_{p'}\|_2 \lambda_{p'}^{-2}\lambda_{q-p'}^{s-1}\\
\leq & c_r \kappa \sum_{p > Q_{b,h}+2}\sum_{|p-q| \leq 2}\lambda_q^{s+1}\|w_q\|_2 \sum_{Q_{b,h}< p' \leq p-2}\lambda_{p'}^{s+1} \|m_{p'}\|_2\lambda_{q-p'}^{s-1}\\
\leq & c_r \kappa \sum_{q \geq -1 }(\lambda_q^{2s+2}\|w_q\|_2^2+\lambda_q^{2s+2}\|m_q\|_2^2).
\end{split}
\end{equation}

One splits $C_3$ into lower and higher modes. 
\begin{equation}\notag
\begin{split}
C_3 = & \sum_{Q_{b,h}-1 \leq p \leq Q_{b,h}} \sum_{q \leq p+2} \lambda_q^{2s}\int_{\T^3}\Delta_q(b_p \cdot \nabla \tilde m_p)w_q \mathrm{d}x\\
& + \sum_{p > Q_{b,h}} \sum_{q \leq p+2} \lambda_q^{2s}\int_{\T^3}\Delta_q(b_p \cdot \nabla \tilde m_p)w_q \mathrm{d}x\\
=:& C_{31}+C_{32}.
\end{split}
\end{equation}

$C_{31},$ made up from the scarce lower modes, is estimated as follows.
\begin{equation}\notag
\begin{split}
|C_{31}| \leq & \sum_{Q_{b,h}-1 \leq p \leq Q_{b,h}} \|b_p\|_\infty \| \nabla \tilde m_p\|_2 \sum_{q \leq p+2}\lambda_q^{2s}\|w_q\|_2\\
\leq & c_r \kappa \sum_{Q_{b,h}-1 \leq p \leq Q_{b,h}}\lambda_p^{s+1}\|m_p\|_2 \sum_{q \leq p+2}\lambda_q^{s+1}\|w_q\|_2\lambda_{q-p}^{s}\\
\leq & c_r \kappa \sum_{q \geq -1 }(\lambda_q^{2s+2}\|w_q\|_2^2+\lambda_q^{2s+2}\|m_q\|_2^2).
\end{split}
\end{equation}

One recalls Definition \ref{wvnbr} and applies H\"older's, Young's and Jensen's inequalities to bound $C_{32}.$
\begin{equation}\notag
\begin{split}
|C_{32}| \leq & \sum_{p > Q_{b,h}}\|b_p\|_\infty \| \nabla \tilde m_p\|_2  \sum_{q \leq p+2}\lambda_q^{2s}\|w_q\|_2\\
\leq & c_r \kappa \sum_{p > Q_{b,h}}\lambda_p^{s+1}\|m_p\|_2 \sum_{q \leq p+2}\lambda_q^{s+1}\|w_q\|_2\lambda_{q-p}^{s}\\
\leq & c_r \kappa \sum_{q \geq -1 }(\lambda_q^{2s+2}\|w_q\|_2^2+\lambda_q^{2s+2}\|m_q\|_2^2).
\end{split}
\end{equation}

\subsection{Estimation of D}

Bony's paraproduct decomposition yields
\begin{equation}\notag
\begin{split}
D= & \sum_{q \geq -1}\sum_{|p-q| \leq 2} \lambda_q^{2s}\int_{\T^3} \Delta_q(m_p \cdot \nabla h_{\leq p-2})w_q \mathrm{d}x\\
&+ \sum_{q \geq -1}\sum_{|p-q| \leq 2} \lambda_q^{2s}\int_{\T^3} \Delta_q(m_{\leq p-2}\cdot \nabla h_p)w_q \mathrm{d}x\\
&+ \sum_{q \geq -1}\sum_{p \geq q-2} \lambda_q^{2s}\int_{\T^3} \Delta_q(\tilde m_p \cdot \nabla h_p)w_q \mathrm{d}x\\
=: & D_1+ D_2 + D_3.
\end{split}
\end{equation}

Utilizing the wavenumber $Q_{b,h},$ one splits $D_1$ into three terms.
\begin{equation}\notag
\begin{split}
D_1=&\sum_{Q_{b,h}< p \leq Q_{b,h}+2}\sum_{|q-p| \leq 2} \lambda_q^{2s}\int_{\T^3} \Delta_q(m_p \cdot \nabla h_{\leq Q_{b,h}})w_q \mathrm{d}x\\
&+\sum_{p> Q_{b,h}+2}\sum_{|q-p| \leq 2} \lambda_q^{2s}\int_{\T^3} \Delta_q(m_p \cdot \nabla h_{\leq Q_{b,h}})w_q \mathrm{d}x\\
&+ \sum_{p> Q_{b,h}+2}\sum_{|q-p| \leq 2} \lambda_q^{2s}\int_{\T^3} \Delta_q(m_p \cdot \nabla h_{(Q_{b,h}, p-2]})w_q \mathrm{d}x\\
=: & D_{11}+D_{12}+D_{13}
\end{split}
\end{equation}

One can estimate $|D_{11}|+|D_{12}|$ without difficulties.
\begin{equation}\notag
\begin{split}
|D_{11}|+|D_{12}| \leq & \| \nabla h_{\leq Q_{b,h}}\|_\infty \sum_{p >Q_{b,h}}\|m_p\|_2 \sum_{|q-p| \leq 2}\lambda_q^{2s}\|w_q\|_2 \\
\leq & c_r\kappa\sum_{p >Q_{b,h}}\lambda_p\|m_p\|_2 \sum_{|q-p| \leq 2}\lambda_q^{2s}\|w_q\|_2 \\
\leq & c_r\kappa\sum_{p >Q_{b,h}}\lambda_p^{s+1}\|m_p\|_2 \sum_{|q-p| \leq 2}\lambda_q^{s+1}\|w_q\|_2 \\
\leq & c_r \kappa \sum_{q \geq -1 }(\lambda_q^{2s+2}\|w_q\|_2^2+\lambda_q^{2s+2}\|m_q\|_2^2).
\end{split}
\end{equation}

By Definition \ref{wvnbr}, H\"older's, Young's and Jensen's inequalities, one has
\begin{equation}\notag
\begin{split}
|D_{13}| \leq & \sum_{p> Q_{b,h}+2}\|m_p\|_2 \sum_{|q-p| \leq 2}\lambda_q^{2s}\|w_q\|_2 \sum_{ Q_{b,h}< p'\leq p-2}\lambda_{p'}\|h_{p'}\|_\infty \\
\leq & c_r \kappa \sum_{p> Q_{b,h}+2}\lambda_{p}^2\|m_p\|_2 \sum_{|q-p| \leq 2}\lambda_q^{2s}\|w_q\|_2 \sum_{ Q_{b,h}< p'\leq p-2}\lambda_{p'-p}^2 \\
\lesssim & c_r \kappa \sum_{p> Q_{b,h}+2}\lambda_{p}^{s+1}\|m_p\|_2 \sum_{|q-p| \leq 2}\lambda_q^{s+1}\|w_q\|_2 \sum_{ Q_{b,h}< p'\leq p-2}\lambda_{p'-p}^2 \\
\lesssim & c_r \kappa \sum_{q \geq -1 }(\lambda_q^{2s+2}\|w_q\|_2^2+\lambda_q^{2s+2}\|m_q\|_2^2).
\end{split}
\end{equation}

It turns out that $D_2$ consists of only higher modes, as $m_{\leq Q_{b,h}}=0.$
\begin{equation}\notag
\begin{split}
D_2 = & \sum_{p > Q_{b,h}+2}\sum_{|q-p| \leq 2} \lambda_q^{2s}\int_{\T^3} \Delta_q(m_{(Q_{b,h}, p-2]}\cdot \nabla h_p)w_q \mathrm{d}x.
\end{split}
\end{equation}

Using H\"older's Young's and Jensen's inequalities, one estimates $D_2.$ 
\begin{equation}\notag
\begin{split}
|D_2| \leq & \sum_{p > Q_{b,h}+2}\lambda_p\|h_p\|_\infty\sum_{|q-p| \leq 2} \lambda_q^{2s}\|w_q\|_2 \sum_{Q_{b,h} < p' \leq p-2}\|m_{p'}\|_2\\
\lesssim & c_r \kappa \sum_{q > Q_{b,h}}\lambda_q^{2s+1}\|w_q\|_2 \sum_{Q_{b,h} < p' \leq q}\|m_{p'}\|_2\\
\lesssim & c_r \kappa \sum_{q > Q_{b,h}}\lambda_q^{s+1}\|w_q\|_2 \sum_{Q_{b,h} < p' \leq q}\lambda_{p'}^{s+1}\|m_{p'}\|_2\lambda_{q-p'}^{s-1}\\
\lesssim & c_r \kappa \sum_{q \geq -1 }(\lambda_q^{2s+2}\|w_q\|_2^2+\lambda_q^{2s+2}\|m_q\|_2^2).
\end{split}
\end{equation}

One split $D_3$ into the lower modes, which are rather few, and the higher modes, which are the majority.
\begin{equation}\notag
\begin{split}
D_3 = &\sum_{q \leq Q_{b,h}+2} \lambda_q^{2s}\int_{\T^3} \Delta_q(m_{Q_{b,h}+1} \cdot \nabla h_{Q_{b,h}})w_q \mathrm{d}x\\
& +\sum_{p > Q_{b,h}}\sum_{q \leq p+2} \lambda_q^{2s}\int_{\T^3} \Delta_q(\tilde m_p \cdot \nabla h_p)w_q \mathrm{d}x\\
=:& D_{31}+D_{32}.
\end{split}
\end{equation}

$D_{31}$ satisfies the following estimate.
\begin{equation}\notag
\begin{split}
D_{31} \leq & \|\nabla h_{Q_{b,h}}\|_\infty \sum_{-1 \leq  q \leq Q_{b,h}+2} \lambda_{q}^{s}\|w_q\|_2\|m_{Q_{b,h}+1}\|_2\\
\lesssim & c_r\kappa \lambda_{Q_{b,h}+1}^{s+1}\|m_{Q_{b,h}+1}\|_2  \sum_{-1 \leq  q \leq Q_{b,h}+2} \lambda_{q}^{2s+1}\|w_q\|_2\\
\lesssim & c_r \kappa \sum_{q \geq -1 }(\lambda_q^{2s+2}\|w_q\|_2^2+\lambda_q^{2s+2}\|m_q\|_2^2).
\end{split}
\end{equation}

The estimate for $D_{32}$ follows from Definition \ref{wvnbr}, H\"older's, Young's and Jensen's inequalities.
\begin{equation}\notag
\begin{split}
D_{32}\leq & \sum_{p > Q_{b,h}}\|\nabla h_p\|_\infty \|\tilde m_p\|_2\sum_{q \leq p+2} \lambda_q^{2s}\|w_q\|_2\\
\lesssim & c_r \kappa \sum_{p > Q_{b,h}}\lambda_p\|\tilde m_p\|_2\sum_{q \leq p+2} \lambda_q^{2s}\|w_q\|_2\\
\lesssim & c_r \kappa \sum_{p > Q_{b,h}}\lambda_p^{s+1}\|m_p\|_2\sum_{q \leq p+2} \lambda_q^{s+1}\|w_q\|_2 \lambda_{q-p}^{s-1}\\
\lesssim & c_r \kappa \sum_{q \geq -1 }(\lambda_q^{2s+2}\|w_q\|_2^2+\lambda_q^{2s+2}\|m_q\|_2^2).
\end{split}
\end{equation}

\subsection{Estimation of E}

One decomposes $E$ using Bony's paraproduct.
\begin{equation}\notag
\begin{split}
E=& -\sum_{q \geq -1}\sum_{|p-q| \leq 2} \lambda_q^{2s}\int_{\T^3}  \Delta_q(v_{p}\cdot \nabla m_{\leq p-2})m_q \mathrm{d}x\\
&-\sum_{q \geq -1}\sum_{|p-q| \leq 2} \lambda_q^{2s}\int_{\T^3}  \Delta_q(v_{\leq p-2}\cdot \nabla m_p)m_q \mathrm{d}x\\
&-\sum_{q \geq -1}\sum_{p \geq q-2} \lambda_q^{2s}\int_{\T^3}  \Delta_q(v_p \cdot \nabla \tilde m_p)m_q \mathrm{d}x\\
=: & E_1 + E_2 +E_3.
\end{split}
\end{equation}

Utilizing the wavenumber $Q_{u,v},$ $E_1$ is split into two.
\begin{equation}\notag
\begin{split}
E_1 =& -\sum_{p \leq Q_{u,v}}\sum_{|q-p| \leq 2} \lambda_q^{2s}\int_{\T^3}  \Delta_q(v_{p}\cdot \nabla m_{\leq p-2})m_q \mathrm{d}x\\
& -\sum_{p > Q_{u,v}}\sum_{|q-p| \leq 2} \lambda_q^{2s}\int_{\T^3}  \Delta_q(v_{p}\cdot \nabla m_{\leq p-2})m_q \mathrm{d}x\\
=:& E_{11}+E_{12}.
\end{split}
\end{equation}

By Definition \ref{uvnbr}, H\"older's, Young's and Jensen's inequalities, $E_{11}$ and $E_{12}$ are estimated in the following ways.
\begin{equation}\notag
\begin{split}
|E_{11}| \leq & \sum_{p \leq Q_{u,v}}\|v_p\|_r \|\nabla m_{\leq p-2} \|_\frac{2r}{r-2}\sum_{|q-p| \leq 2}\lambda_q^{2s}\|m_q\|_2\\
\leq & \sum_{p \leq Q_{u,v}}\lambda_p^{-1+\frac{3}{r}}\|v_p\|_r\sum_{p' \leq p-2}\lambda_{p'}^{1+\frac{3}{r}}\|m_{p'} \|_{2}\sum_{|q-p| \leq 2}\lambda_q^{2s+1-\frac{3}{r}}\|m_q\|_2\\
\lesssim & c_r\kappa \sum_{q \leq Q_{u,v}+2}\lambda_q^{s+1}\|m_q\|_2\sum_{p' \leq q}\lambda_{p'}^{s+1}\|m_{p'} \|_{2}\lambda_{q-p'}^{s-\frac{3}{r}}\\
\lesssim & c_r\kappa \sum_{q \geq -1}\lambda_q^{2s+2} \|m_q\|_2^2;
\end{split}
\end{equation}
\begin{equation}\notag
\begin{split}
|E_{12}| \leq & \sum_{p > Q_{u,v}}\|v_p\|_r \|\nabla m_{\leq p-2} \|_{\frac{2r}{r-2}}\sum_{|q-p| \leq 2}\lambda_q^{2s}\|m_q\|_2\\
\leq & \sum_{p > Q_{u,v}}\lambda_p^{-1+\frac{3}{r}}\|v_p\|_r\sum_{p' \leq p-2}\lambda_{p'}^{1+\frac{3}{r}}\|m_{p'} \|_{2}\sum_{|q-p| \leq 2}\lambda_q^{2s+1-\frac{3}{r}}\|m_q\|_2\\
\lesssim & c_r\kappa \sum_{q > Q_{u,v}-2}\lambda_q^{s+1}\|m_q\|_2\sum_{p' \leq q}\lambda_{p'}^{s+1}\|m_{p'} \|_{2}\lambda_{q-p'}^{s-\frac{3}{r}}\\
\lesssim & c_r\kappa \sum_{q \geq -1}\lambda_q^{2s+2} \|m_q\|_2^2.
\end{split}
\end{equation}

By the commutator notation, one has
\begin{equation}\notag
\begin{split}
E_{2}= & \sum_{q \geq -1} \sum_{|p-q| \leq 2} \lambda_q^{2s}\int_{\T^3}[\Delta_q, v_{\leq p-2}\cdot \nabla] m_p m_q \mathrm{d}x\\
&+ \sum_{q \geq -1} \sum_{|p-q| \leq 2} \lambda_q^{2s}\int_{\T^3}v_{\leq q-2} \cdot \nabla \Delta_q m_p m_q \mathrm{d}x\\
&+ \sum_{q \geq -1}\sum_{|p-q| \leq 2} \lambda_q^{2s}\int_{\T^3}(v_{\leq p-2}-v_{\leq q-2})\cdot \nabla \Delta_q m_p m_q \mathrm{d}x\\
=: & E_{21}+E_{22}+E_{23},
\end{split}
\end{equation}
where $E_{22}$ vanishes as $\nabla \cdot v_{\leq q-2}=0.$

Splitting $E_{21}$ by the wavenumber $Q_{u,v},$ one has
\begin{equation}\notag
\begin{split}
E_{21} = & \sum_{-1 \leq p \leq Q_{u,v}+2} \sum_{|q-p| \leq 2} \lambda_q^{2s}\int_{\T^3}[\Delta_q, v_{\leq p-2}\cdot \nabla] m_p m_q \mathrm{d}x\\
&+ \sum_{p > Q_{u,v}+2} \sum_{|q-p| \leq 2} \lambda_q^{2s}\int_{\T^3}[\Delta_q, v_{\leq Q_{u,v}}\cdot \nabla] m_p m_q \mathrm{d}x\\
&+ \sum_{p > Q_{u,v}+2} \sum_{|q-p| \leq 2} \lambda_q^{2s}\int_{\T^3}[\Delta_q, v_{(Q_{u,v}, p-2]}\cdot \nabla] m_p m_q \mathrm{d}x\\
=: & E_{211}+E_{212}+E_{213}.
\end{split}
\end{equation}

Using Definition \ref{uvnbr}, Lemma \ref{cmest}, H\"older's and Young's inequalities, one can estimate $E_{211}.$
\begin{equation}\notag
\begin{split}
|E_{211}| \leq & \sum_{-1 \leq p \leq Q_{u,v}+2}\|m_p\|_\frac{2r}{r-2} \sum_{|q-p| \leq 2} \lambda_q^{2s} \|m_q\|_2 \sum_{p' \leq p-2}\lambda_{p'}\|v_{p'}\|_r\\
\lesssim & c_r \kappa \sum_{-1 \leq p \leq Q_{u,v}+2}\lambda_p^{\frac{3}{r}}\|m_p\|_2 \sum_{|q-p| \leq 2} \lambda_q^{2s} \|m_q\|_2 \lambda_{p}^{2-\frac{3}{r}}\\
\lesssim & c_r \kappa \sum_{-1 \leq p \leq Q_{u,v}+2}\lambda_p^{s+1}\|m_p\|_2 \sum_{|q-p| \leq 2} \lambda_q^{s+1} \|m_q\|_2 \\
\lesssim & c_r\kappa \sum_{q \geq -1}\lambda_q^{2s+2} \|m_q\|_2^2.
\end{split}
\end{equation}

The term $E_{212}$ can be estimated in a similar fashion.
\begin{equation}\notag
\begin{split}
|E_{212}| \leq & \sum_{p > Q_{u,v}+2}\|m_p\|_\frac{2r}{r-2} \sum_{|q-p| \leq 2} \lambda_q^{2s} \|m_q\|_2 \sum_{p' \leq Q_{u,v}}\lambda_{p'}\|v_{p'}\|_r\\
\lesssim & c_r \kappa \sum_{p > Q_{u,v}+2}\lambda_p^{\frac{3}{r}}\|m_p\|_2 \sum_{|q-p| \leq 2} \lambda_q^{2s} \|m_q\|_2 \lambda_{p}^{2-\frac{3}{r}}\\
\lesssim & c_r \kappa \sum_{p > Q_{u,v}+2}\lambda_p^{s+1}\|m_p\|_2 \sum_{|q-p| \leq 2} \lambda_q^{s+1} \|m_q\|_2 \\
\lesssim & c_r\kappa \sum_{q \geq -1}\lambda_q^{2s+2} \|m_q\|_2^2.
\end{split}
\end{equation}

The estimate for $E_{213}$ follows from Definition \ref{uvnbr}, Lemma \ref{cmest}, H\"older's and Young's inequalities.
\begin{equation}\notag
\begin{split}
|E_{213}| \leq & \sum_{p > Q_{u,v}+2}\|m_p\|_\frac{2r}{r-2} \sum_{|q-p| \leq 2} \lambda_q^{2s} \|m_q\|_2 \sum_{Q_{u,v} < p' \leq p-2}\lambda_{p'}\|v_{p'}\|_r\\
\leq & c_r \kappa \sum_{p > Q_{u,v}+2}\lambda_p^{\frac{3}{r}}\|m_p\|_2 \sum_{|q-p| \leq 2} \lambda_q^{2s} \|m_q\|_2 \sum_{Q_{u,v} < p' \leq p-2}\lambda_{p'}^{2-\frac{3}{r}}\\
\lesssim & c_r \kappa \sum_{p > Q_{u,v}+2}\lambda_p^{s+1}\|m_p\|_2 \sum_{|q-p| \leq 2} \lambda_q^{s+1} \|m_q\|_2 \sum_{Q_{u,v} < p' \leq p-2}\lambda_{p'-p}^{2-\frac{3}{r}}\\
\lesssim & c_r\kappa \sum_{q \geq -1}\lambda_q^{2s+2} \|m_q\|_2^2.
\end{split}
\end{equation}

Explicitly writing out $(v_{\leq p-2}- v_{\leq q-2})$ leads to
\begin{equation}\notag
\begin{split}
|E_{23}| \lesssim & \sum_{q \geq -1}\sum_{|p-q| \leq 2} \lambda_q^{2s}\int_{\T^3}(|v_{q-3}|+|v_{q-2}|+ |v_{q-1}|+ |v_{q}|)|\nabla \Delta_q m_p m_q| \mathrm{d}x\\
\lesssim & \sum_{-1 \leq q \leq Q_{u,v}}\sum_{|p-q| \leq 2} \lambda_q^{2s}\int_{\T^3}|v_{q}||\nabla \Delta_q m_p m_q| \mathrm{d}x\\
& + \sum_{q > Q_{u,v}}\sum_{|p-q| \leq 2} \lambda_q^{2s}\int_{\T^3}|v_{q}||\nabla \Delta_q m_p m_q| \mathrm{d}x\\
=:& E_{231}+E_{232}.
\end{split}
\end{equation}

The estimate for $E_{231}$ is as follows.
\begin{equation}\notag
\begin{split}
E_{231} \lesssim & \sum_{-1 \leq q \leq Q_{u,v}}\lambda_q^{2s}\|v_q\|_r \|m_q\|_\frac{2r}{r-2}\sum_{|p-q| \leq 2} \|\nabla m_p\|_2\\
\lesssim & c_r \kappa \sum_{-1 \leq q \leq Q_{u,v}}\lambda_q^{2s+1} \|m_q\|_2\sum_{|p-q| \leq 2}\lambda_p\|m_p\|_2\\
\lesssim & c_r \kappa \sum_{-1 \leq q \leq Q_{u,v}}\lambda_q^{s+1}\|m_q\|_2 \sum_{|p-q| \leq 2} \lambda_p^{s+1} \|m_p\|_2\\
\lesssim & c_r\kappa \sum_{q \geq -1}\lambda_q^{2s+2} \|m_q\|_2^2.
\end{split}
\end{equation}

By Definition \ref{uvnbr}, H\"older's and Young's inequalities, $E_{232}$ can be estimated.
\begin{equation}\notag
\begin{split}
E_{232} \leq & \sum_{q > Q_{u,v}}\lambda_q^{2s}\|v_q\|_r \|m_q\|_\frac{2r}{r-2}\sum_{|p-q| \leq 2} \|\nabla m_p\|_2\\
\lesssim & c_r \kappa \sum_{q > Q_{u,v}}\lambda_q^{2s+1} \|m_q\|_2\sum_{|p-q| \leq 2}\lambda_p\|m_p\|_2\\
\lesssim & c_r \kappa \sum_{q > Q_{u,v}}\lambda_q^{s+1}\|m_q\|_2 \sum_{|p-q| \leq 2} \lambda_p^{s+1} \|m_p\|_2\\
\lesssim & c_r\kappa \sum_{q \geq -1}\lambda_q^{2s+2} \|m_q\|_2^2.
\end{split}
\end{equation}

One separates the lower and higher modes of $E_3$ using the wavenumber $Q_{u,v}.$
\begin{equation}\notag
\begin{split}
E_{3}=& -\sum_{-1 \leq p \leq Q_{u,v}}\sum_{q \leq p-2} \lambda_q^{2s}\int_{\T^3}  \Delta_q(v_p \cdot \nabla \tilde m_p)m_q \mathrm{d}x\\
& -\sum_{p > Q_{u,v}}\sum_{q \leq p-2} \lambda_q^{2s}\int_{\T^3}  \Delta_q(v_p \cdot \nabla \tilde m_p)m_q \mathrm{d}x\\
=: & E_{31} + E_{32}.
\end{split}
\end{equation}

With the help of Definition \ref{uvnbr}, H\"older's, Young's and Jensen's inequalities, the terms $E_{31}$ and $E_{32}$ can be under control.
\begin{equation}\notag
\begin{split}
|E_{31}| \leq &  \sum_{-1 \leq p \leq Q_{u,v}}\|v_p\|_r\|\nabla m_p\|_2\sum_{q \leq p-2}\lambda_q^{2s}\|m_q\|_\frac{2r}{r-2}\\
\leq &  c_r\kappa \sum_{-1 \leq p \leq Q_{u,v}}\lambda_p^{2-\frac{3}{r}}\|m_p\|_2\sum_{q \leq p-2}\lambda_q^{2s}\|m_q\|_\frac{2r}{r-2}\\
\leq &  c_r\kappa \sum_{-1 \leq p \leq Q_{u,v}}\lambda_p^{2-\frac{3}{r}}\|m_p\|_2\sum_{q \leq p-2}\lambda_q^{2s+\frac{3}{r}}\|m_q\|_2\\
\leq &  c_r\kappa \sum_{-1 \leq p \leq Q_{u,v}}\lambda_p^{s+1}\|m_p\|_2\sum_{q \leq p-2}\lambda_q^{s+1}\|m_q\|_2\lambda_{q-p}^{s+\frac{3}{r}-1}\\
\leq & c_r\kappa \sum_{q \geq -1}\lambda_q^{2s+2} \|m_q\|_2^2;
\end{split}
\end{equation}
\begin{equation}\notag
\begin{split}
|E_{32}| \leq &  \sum_{p > Q_{u,v}}\|v_p\|_r\|\nabla m_p\|_2\sum_{q \leq p-2}\lambda_q^{2s}\|m_q\|_\frac{2r}{r-2}\\
\leq &  c_r\kappa \sum_{p > Q_{u,v}}\lambda_p^{2-\frac{3}{r}}\|m_p\|_2\sum_{q \leq p-2}\lambda_q^{2s}\|m_q\|_\frac{2r}{r-2}\\
\leq &  c_r\kappa \sum_{p > Q_{u,v}}\lambda_p^{2-\frac{3}{r}}\|m_p\|_2\sum_{q \leq p-2}\lambda_q^{2s+\frac{3}{r}}\|m_q\|_2\\
\leq &  c_r\kappa \sum_{p > Q_{u,v}}\lambda_p^{s+1}\|m_p\|_2\sum_{q \leq p-2}\lambda_q^{s+1}\|m_q\|_2\lambda_{q-p}^{s+\frac{3}{r}-1}\\
\leq & c_r\kappa \sum_{q \geq -1}\lambda_q^{2s+2} \|m_q\|_2^2.
\end{split}
\end{equation}

\subsection{Estimation of F}

By Bony's paraproduct decomposition, one has
\begin{equation}\notag
\begin{split}
F=& -\sum_{q \geq -1}\sum_{|p-q| \leq 2} \lambda_q^{2s}\int_{\T^3}  \Delta_q(w_{p}\cdot \nabla b_{\leq p-2})m_q \mathrm{d}x\\
&-\sum_{q \geq -1}\sum_{|p-q| \leq 2} \lambda_q^{2s}\int_{\T^3}  \Delta_q(w_{\leq p-2}\cdot \nabla b_p)m_q \mathrm{d}x\\
&-\sum_{q \geq -1}\sum_{p \geq q-2} \lambda_q^{2s}\int_{\T^3}  \Delta_q(\tilde w_p \cdot \nabla b_p)m_q \mathrm{d}x\\
=:& F_1 + F_2 + F_3.
\end{split}
\end{equation}

Using the fact that $m_{\leq Q_{b,h}}=0,$ one splits $F_1$ into two terms.
\begin{equation}\notag
\begin{split}
F_1= & -\sum_{p > Q_{b,h}+2}\sum_{|q-p| \leq 2} \lambda_q^{2s}\int_{\T^3}  \Delta_q(w_{p}\cdot \nabla b_{\leq Q_{b,h}})m_q \mathrm{d}x\\
& -\sum_{p > Q_{b,h}+2}\sum_{|q-p| \leq 2} \lambda_q^{2s}\int_{\T^3}  \Delta_q(w_{p}\cdot \nabla b_{(Q_{b,h}, p-2]})m_q \mathrm{d}x\\
=:& F_{11}+F_{12}.
\end{split}
\end{equation}

To estimate $F_{11},$ one uses Definition \ref{wvnbr}, H\"older's and Young's inequalities.
\begin{equation}\notag
\begin{split}
|F_{11}| \leq & \|\nabla b_{\leq Q_{b,h}}\|_\infty \sum_{q > Q_{b,h}}\lambda_q^{2s}\|m_q\|_2\sum_{|p-q| \leq 2}\|w_p\|_2\\
\leq & \|b_{\leq Q_{b,h}}\|_\infty \sum_{q > Q_{b,h}}\lambda_q^{s+1}\|m_q\|_2\sum_{|p-q| \leq 2}\lambda_p^{s}\|w_p\|_2\\
\leq & c_r \kappa \sum_{q > Q_{b,h}}\lambda_q^{s+1}\|m_q\|_2\sum_{|p-q| \leq 2}\lambda_p^{s+1}\|w_p\|_2\\
\leq & c_r \kappa \sum_{q > -1 }(\lambda_q^{2s+2}\|w_q\|_2^2+\lambda_q^{2s+2}\|m_q\|_2^2).
\end{split}
\end{equation}

By Definition \ref{wvnbr}, H\"older's and Young's inequalities, $F_{12}$ satisfies the following.
\begin{equation}\notag
\begin{split}
|F_{12}| \leq & \sum_{q > Q_{b,h}}\lambda_q^{2s}\|m_q\|_2\sum_{|p-q| \leq 2}\|w_p\|_2\sum_{Q_{b,h}<p' \leq p-2}\lambda_{p'}\|b_{p'}\|_\infty\\
\leq & \sum_{q > Q_{b,h}}\lambda_q^{2s+1}\|m_q\|_2\sum_{|p-q| \leq 2}\lambda_p\|w_p\|_2\sum_{Q_{b,h}<p' \leq p-2}\lambda_{p'-p}\lambda_{p}^{-1}\|b_{p'}\|_\infty\\
\leq & c_r \kappa \sum_{q > Q_{b,h}}\lambda_q^{s+1}\|m_q\|_2\sum_{|p-q| \leq 2}\lambda_p^{s+1}\|w_p\|_2\sum_{Q_{b,h}<p' \leq p-2}\lambda_{p'-p}^2\\
\leq & c_r \kappa \sum_{q > -1 }(\lambda_q^{2s+2}\|w_q\|_2^2+\lambda_q^{2s+2}\|m_q\|_2^2).
\end{split}
\end{equation}

$F_2$ is split into lower and higher modes based on the wavenumber $Q_{b,h}$ as well as the fact that $m_{\leq Q_{b,h}}=0.$
\begin{equation}\notag
\begin{split}
F_2= & -\sum_{Q_{b,h}-2< p \leq Q_{b,h}}\sum_{|q-p| \leq 2} \lambda_q^{2s}\int_{\T^3}  \Delta_q(w_{\leq p-2}\cdot \nabla b_p)m_q \mathrm{d}x\\
& -\sum_{p > Q_{b,h}}\sum_{|q-p| \leq 2} \lambda_q^{2s}\int_{\T^3}  \Delta_q(w_{\leq p-2}\cdot \nabla b_p)m_q \mathrm{d}x\\
=: & F_{21}+F_{22}.
\end{split}
\end{equation}

It follows from Definition \ref{wvnbr}, H\"older's, Young's and Jensen's inequalities that
\begin{equation}\notag
\begin{split}
|F_{21}| \leq & \sum_{Q_{b,h}-2< p \leq Q_{b,h}}\lambda_p\|b_p\|_\infty \sum_{|q-p| \leq 2} \lambda_q^{2s}\|m_q\|_2\sum_{p' \leq p-2}\|w_{p'}\|_2\\
\lesssim & c_r \kappa \sum_{Q_{b,h}< q \leq Q_{b,h}+2} \lambda_q^{2s+1}\|m_q\|_2\sum_{p' \leq q}\|w_{p'}\|_2\\
\lesssim & c_r \kappa \sum_{Q_{b,h}< q \leq Q_{b,h}+2}\lambda_q^{s+1}\|m_q\|_2\sum_{p' \leq q}\lambda_{p'}^{s+1}\|w_{p'}\|_2\lambda_{p'}^{-s-1}\lambda_q^{s}\\
\leq & c_r \kappa \sum_{q > -1 }(\lambda_q^{2s+2}\|w_q\|_2^2+\lambda_q^{2s+2}\|m_q\|_2^2).
\end{split}
\end{equation}

One estimates $F_{22}$ with the help of H\"older's, Young's and Jensen's inequalities.
\begin{equation}\notag
\begin{split}
|F_{22}| \leq & \sum_{p > Q_{b,h}}\|\nabla b_p\|_\infty \sum_{|q-p| \leq 2}\lambda_q^{2s} \|m_q\|_2\sum_{p' \leq p-2}\|w_{p'}\|_2\\
\leq & \sum_{q > Q_{b,h}}\lambda_q^{2s} \|m_q\|_2\lambda^{1-\delta}_p\Lambda_{b,h}^{\delta}\sum_{|p-q| \leq 2}\lambda^\delta_{p-Q_{b,h}}\|b_p\|_\infty\sum_{ p' \leq p-2} \lambda_{p'}\|w_{p'}\|_2 \lambda_{p'}^{-1}\\
\leq & c_r \kappa \sum_{q > Q_{b,h}}\lambda_q^{s+1} \|m_q\|_2\sum_{ -1 \leq p' \leq q} \lambda_{p'}^{s+1}\|w_{p'}\|_2\lambda_{q-p'}^{s}\lambda_{p'}^{-1}\\
\leq & c_r \kappa \sum_{q > -1 }(\lambda_q^{2s}\|w_q\|_2^2+\lambda_q^{2s+2}\|m_q\|_2^2).
\end{split}
\end{equation}

As $m_{\leq Q_{b,h}}=0,$ one splits $F_3$ into two terms.
\begin{equation}\notag
\begin{split}
F_3= &-\sum_{p \leq Q_{b,h}}\sum_{Q_{b,h}< q \leq p+2} \lambda_q^{2s}\int_{\T^3}  \Delta_q(\tilde w_p \cdot \nabla b_p)m_q \mathrm{d}x\\
&-\sum_{p > Q_{b,h}}\sum_{Q_{b,h}< q \leq p+2} \lambda_q^{2s}\int_{\T^3}  \Delta_q(\tilde w_p \cdot \nabla b_p)m_q \mathrm{d}x\\
=:& F_{31}+F_{32}.
\end{split}
\end{equation}

The estimate for $F_{31}$ is as follows.
\begin{equation}\notag
\begin{split}
|F_{31}| \leq& \sum_{p \leq Q_{b,h}}\|\nabla b_p\|_\infty \|\tilde w_p\|_2\sum_{Q_{b,h}< q \leq p+2} \lambda_q^{2s} \|m_q\|_2\\
\leq& c_r \kappa \sum_{p \leq Q_{b,h}}\lambda_p \|w_p\|_2\sum_{Q_{b,h}< q \leq p+2} \lambda_q^{2s} \|m_q\|_2\\
\leq& c_r \kappa \sum_{p \leq Q_{b,h}}\lambda_p^{s+1} \|w_p\|_2\sum_{Q_{b,h}< q \leq p+2} \lambda_q^{s+1} \|m_q\|_2\lambda_{q-p}^{s}\lambda_q^{-1}\\
\leq & c_r \kappa \sum_{q > -1 }(\lambda_q^{2s+2}\|w_q\|_2^2+\lambda_q^{2s+2}\|m_q\|_2^2).
\end{split}
\end{equation}

One uses H\"older's, Young's and Jensen's inequalities to estimate $F_{32}.$
\begin{equation}\notag
\begin{split}
|F_{32}| \leq& \sum_{p > Q_{b,h}}\|\nabla b_p\|_\infty \|\tilde w_p\|_2\sum_{q \leq p+2} \lambda_q^{2s} \|m_q\|_2\\
\leq& c_r\kappa \sum_{p > Q_{b,h}}\lambda_p\|w_p\|_2\sum_{q \leq p+2} \lambda_q^{2s} \|m_q\|_2\\
\leq& c_r\kappa \sum_{p > Q_{b,h}}\lambda_p^{s+1}\|w_p\|_2\sum_{q \leq p+2} \lambda_q^{s+1} \|m_q\|_2\lambda_{q-p}^{s}\lambda_q^{-1}\\
\leq & c_r \kappa \sum_{q > -1 }(\lambda_q^{2s+2}\|w_q\|_2^2+\lambda_q^{2s+2}\|m_q\|_2^2).
\end{split}
\end{equation}

\subsection{Estimation of G}

Using Bony's paraproduct decomposition, one has
\begin{equation}\notag
\begin{split}
G=& \sum_{q \geq -1}\sum_{|p-q| \leq 2} \lambda_q^{2s}\int_{\T^3}  \Delta_q(b_{p}\cdot \nabla w_{\leq p-2})m_q \mathrm{d}x\\
&+\sum_{q \geq -1}\sum_{|p-q| \leq 2} \lambda_q^{2s}\int_{\T^3}  \Delta_q(b_{\leq p-2}\cdot \nabla w_p)m_q \mathrm{d}x\\
&+\sum_{q \geq -1}\sum_{p \geq q-2} \lambda_q^{2s}\int_{\T^3}  \Delta_q(b_p \cdot \nabla \tilde w_p)m_q \mathrm{d}x\\
=: & G_1+G_2+G_3.
\end{split}
\end{equation}

Taking into account that $m_{\leq Q_{b,h}} =0,$ one separates lower and higher modes of $G_1$ by the wavenumber $Q_{b,h}.$
\begin{equation}\notag
\begin{split}
G_1=& \sum_{Q_{b,h}-2 \leq p \leq Q_{b,h}}\sum_{|q-p| \leq 2} \lambda_q^{2s}\int_{\T^3}  \Delta_q(b_{p}\cdot \nabla w_{\leq p-2})m_q \mathrm{d}x\\
&+\sum_{p > Q_{b,h}}\sum_{|q-p| \leq 2} \lambda_q^{2s}\int_{\T^3}  \Delta_q(b_{p}\cdot \nabla w_{\leq p-2})m_q \mathrm{d}x\\
=:& G_{11}+G_{12}.
\end{split}
\end{equation}

Thanks to the fact that $q=Q_{b,h}+1$ or $Q_{b,h}+2$, one can control $G_{11}.$
\begin{equation}\notag
\begin{split}
|G_{11}| \leq & \sum_{Q_{b,h}-2 < p \leq Q_{b,h}}\|b_p\|_\infty \sum_{-1 \leq p' \leq p-2}\lambda_{p'}\|w_{p'}\|_2\sum_{Q_{b,h}< q \leq Q_{b,h}+2}\lambda_q^{2s}\|m_q\|_2\\
\leq & c_r \kappa \sum_{Q_{b,h}< q \leq Q_{b,h}+2}\lambda_q^{s+1}\|m_q\|_2\lambda_{q}^{s-1}\sum_{-1 \leq p' \leq q}\lambda_{p'}^{s+1}\|w_{p'}\|_2\lambda_{p'}^{-s}\\
\lesssim & c_r \kappa \sum_{q \geq -1 }(\lambda_q^{2s+2}\|w_q\|_2^2+\lambda_q^{2s+2}\|m_q\|_2^2).
\end{split}
\end{equation}

Using Definition \ref{wvnbr}, H\"older's, Young's and Jensen's inequalities, one has
\begin{equation}\notag
\begin{split}
|G_{12}| \leq & \sum_{p > Q_{b,h}} \|b_p\|_\infty \sum_{-1 \leq p' \leq p-2}\lambda_{p'}\|w_{p'}\|_2\sum_{|q-p| \leq 2}\lambda_q^{2s}\|m_q\|_2\\
\leq & c_r \kappa \sum_{p > Q_{b,h}}\sum_{-1 \leq p' \leq p-2}\lambda_{p'}^{s+1}\|w_{p'}\|_2\lambda_{p'}^{-s}\sum_{|q-p| \leq 2}\lambda_q^{s+1}\|m_q\|_2\lambda_{q}^{s-1}\\
\lesssim & c_r \kappa \sum_{p > Q_{b,h}}\lambda_p^{s+1}\|m_p\|_2\sum_{-1 \leq p' \leq p-2}\lambda_{p'}^{s+1}\|w_{p'}\|_2\lambda_{p'}^{-s}\lambda_{p}^{s-1}\\
\lesssim & c_r \kappa \sum_{q \geq -1 }(\lambda_q^{2s+2}\|w_q\|_2^2+\lambda_q^{2s+2}\|m_q\|_2^2).
\end{split}
\end{equation}

Rewriting $G_2$ using the commutator notation yields
\begin{equation}\notag
\begin{split}
G_2=& \sum_{q \geq -1}\sum_{|p-q| \leq 2} \lambda_q^{2s}\int_{\T^3}  [\Delta_q, b_{\leq p-2}\cdot \nabla] w_p m_q \mathrm{d}x\\
&+\sum_{q \geq -1 }\sum_{|p-q| \leq 2} \lambda_q^{2s}\int_{\T^3}  b_{ \leq q-2}\cdot \nabla \Delta_q w_p m_q \mathrm{d}x\\
&+\sum_{q \geq -1 }\sum_{|p-q| \leq 2} \lambda_q^{2s}\int_{\T^3}  (b_{p-2}-b_{q-2})\cdot \nabla \Delta_q w_p m_q \mathrm{d}x\\
=:& G_{21}+G_{22}+G_{23}.
\end{split}
\end{equation}

One further splits $G_{21}$ into three parts by the wavenumber $Q_{b,h}.$
\begin{equation}\notag
\begin{split}
G_{21} = &  \sum_{Q_{b,h}-2 < p \leq Q_{b,h}+2}\sum_{|q-p| \leq 2} \lambda_q^{2s}\int_{\T^3}  [\Delta_q, b_{\leq p-2}\cdot \nabla] w_p m_q \mathrm{d}x\\
&+ \sum_{p > Q_{b,h}+2}\sum_{|q-p| \leq 2} \lambda_q^{2s}\int_{\T^3}  [\Delta_q, b_{\leq Q_{b,h}}\cdot \nabla] w_p m_q \mathrm{d}x\\
&+ \sum_{p > Q_{b,h}+2}\sum_{|q-p| \leq 2} \lambda_q^{2s}\int_{\T^3}  [\Delta_q, b_{(Q_{b,h}, p-2]}\cdot \nabla] w_p m_q \mathrm{d}x\\
=:& G_{211}+G_{212}+G_{213}.
\end{split}
\end{equation}

Using Definition \ref{wvnbr}, H\"older's and Young's inequalities, one can estimate $|G_{211}|+|G_{212}|.$
\begin{equation}\notag
\begin{split}
|G_{211}|+|G_{212}| \leq & \|\nabla b_{\leq Q_{b,h}}\|_\infty\sum_{p > Q_{b,h}-2} \|w_p\|_2 \sum_{|q-p| \leq 2}\lambda_q^{2s}\|m_q\|_2\\
\lesssim & \|b_{\leq Q_{b,h}}\|_\infty\sum_{p > Q_{b,h}-2}\lambda_p \|w_p\|_2 \sum_{|q-p| \leq 2}\lambda_q^{2s}\|m_q\|_2\\
\lesssim & c_r \kappa \sum_{p \geq -1}\lambda_p^{s+1} \|w_p\|_2 \sum_{|q-p| \leq 2}\lambda_q^{s+1}\|m_q\|_2\\
\lesssim & c_r \kappa \sum_{p \geq -1 }(\lambda_q^{2s+2}\|w_q\|_2^2+\lambda_q^{2s+2}\|m_q\|_2^2).
\end{split}
\end{equation}

The estimate for $G_{213}$ is as follows.
\begin{equation}\notag
\begin{split}
|G_{213}| \leq & \sum_{p > Q_{b,h}+2}\|w_p\|_2\sum_{Q_{h,b}< p' \leq p-2}\lambda_{p'}\|b_{p'}\|_\infty \sum_{|q-p| \leq 2} \lambda_q^{2s}\| m_q \|_2 \\
\lesssim & c_r \kappa \sum_{q \geq -1}\lambda_p^{s+1} \|w_p\|_2\sum_{|q-p| \leq 2}\lambda_q^{s+1}\|m_q\|_2\lambda_q^{-1}\sum_{Q_{h,b}< p' \leq p-2}\lambda_{p'-p}\\
\lesssim & c_r \kappa \sum_{q \geq -1 }(\lambda_q^{2s+2}\|w_q\|_2^2+\lambda_q^{2s+2}\|m_q\|_2^2).
\end{split}
\end{equation}

As noted before, $G_{22}$ and $C_{12}$ cancel each other.   
\begin{equation}\notag
\begin{split}
C_{12}+G_{22}= & \sum_{q \geq -1 }\sum_{|p-q| \leq 2} \lambda_q^{2s}\int_{\T^3}  b_{ \leq q-2}\cdot \nabla ( \Delta_q w_p m_q + \Delta_q m_p w_q ) \mathrm{d}x\\
= & \sum_{q \geq -1 }\sum_{|p-q| \leq 2} \lambda_q^{2s}\int_{\T^3}  b_{ \leq q-2}\cdot \nabla \Delta_q w_p (m_q+w_q) \mathrm{d}x\\
&+ \sum_{q \geq -1 }\sum_{|p-q| \leq 2} \lambda_q^{2s}\int_{\T^3}  b_{ \leq q-2}\cdot \nabla \Delta_q m_p (w_q+m_q) \mathrm{d}x\\
= & \sum_{q \geq -1 }\lambda_q^{2s}\int_{\T^3}  b_{ \leq q-2}\cdot \nabla (m_q+w_q) (m_q+w_q) \mathrm{d}x\\
=& 0.
\end{split}
\end{equation}

Since $m_{\leq Q_{b,h}}=0,$ $G_{23}$ consists of mostly higher modes.
\begin{equation}\notag
\begin{split}
|G_{23}| \lesssim & \sum_{q \geq -1 }\sum_{|p-q| \leq 2} \lambda_q^{2s}\int_{\T^3}  |b_p||\nabla \Delta_q w_p m_q|\mathrm{d}x\\
\lesssim & \sum_{ Q_{b,h}-2 < p \leq Q_{b,h} }\sum_{|q-p| \leq 2} \lambda_q^{2s}\int_{\T^3}  |b_p||\nabla \Delta_q w_p m_q|\mathrm{d}x\\
&+ \sum_{p > Q_{b,h} }\sum_{|q-p| \leq 2} \lambda_q^{2s}\int_{\T^3}  |b_p||\nabla \Delta_q w_p m_q|\mathrm{d}x\\
=:& G_{231}+G_{232}.
\end{split}
\end{equation}

By Definition \ref{wvnbr}, H\"older's and Young's inequalities, one has 
\begin{equation}\notag
\begin{split}
G_{231} \lesssim & \sum_{ -1 \leq p > Q_{b,h} }\|b_p\|_\infty \lambda_p \|w_p\|_2 \sum_{|q-p| \leq 2}\lambda_q^{2s}\|m_q\|_2\\
\lesssim & c_r \kappa \sum_{p \geq -1}\lambda_p^{s+1} \|w_p\|_2 \sum_{|q-p| \leq 2}\lambda_q^{s+1}\|m_q\|_2\\
\lesssim & c_r \kappa \sum_{q \geq -1 }(\lambda_q^{2s+2}\|w_q\|_2^2+\lambda_q^{2s+2}\|m_q\|_2^2).
\end{split}
\end{equation}

$G_{232}$ is estimated as follows.
\begin{equation}\notag
\begin{split}
G_{232} \lesssim & \sum_{p > Q_{b,h} }\|b_p\|_\infty \lambda_p \|w_p\|_2 \sum_{|q-p| \leq 2}\lambda_q^{2s}\|m_q\|_2\\
\lesssim & c_r \kappa \sum_{p \geq -1}\lambda_p^{s+1} \|w_p\|_2 \sum_{|q-p| \leq 2}\lambda_q^{s+1}\|m_q\|_2\\
\lesssim & c_r \kappa \sum_{q > -1 }(\lambda_q^{2s+2}\|w_q\|_2^2+\lambda_q^{2s+2}\|m_q\|_2^2).
\end{split}
\end{equation}

One divides $G_3$ into lower and higher modes using the wavenumber $Q_{b,h}.$
\begin{equation}\notag
\begin{split}
G_3=& \sum_{Q_{b,h}-2 < p \leq Q_{b,h}}\sum_{q \leq p+2} \lambda_q^{2s}\int_{\T^3}  \Delta_q(b_p \cdot \nabla \tilde w_p)m_q \mathrm{d}x\\
&+\sum_{p > Q_{b,h}}\sum_{q \leq p+2} \lambda_q^{2s}\int_{\T^3}  \Delta_q(b_p \cdot \nabla \tilde w_p)m_q \mathrm{d}x\\
=:& G_{31}+G_{32}.
\end{split}
\end{equation}

One can estimate $G_{31}$ in the following way.
\begin{equation}\notag
\begin{split}
|G_{31}| \leq & \sum_{Q_{b,h}-2< p \leq Q_{b,h}}\|b_p\|_\infty \|\nabla \tilde w_p\|_2 \sum_{q \leq p+2}\lambda_q^{2s}\|m_q\|_2 \\
\lesssim & c_r \kappa \sum_{Q_{b,h}-2< p \leq Q_{b,h}}\lambda_p^{s+1} \|w_p\|_2 \sum_{q \leq p+2}\lambda_q^{s+1}\|m_q\|_2\lambda_{q-p}^{s}\lambda_q^{-1} \\
\lesssim & c_r \kappa \sum_{q \geq -1 }(\lambda_q^{2s+2}\|w_q\|_2^2+\lambda_q^{2s+2}\|m_q\|_2^2).
\end{split}
\end{equation}

Meanwhile, by Definition \ref{wvnbr}, H\"older's, Young's and Jensen's inequalities, it holds that
\begin{equation}\notag
\begin{split}
|G_{32}| \leq & \sum_{p > Q_{b,h}}\|b_p\|_\infty \|\nabla \tilde w_p\|_2 \sum_{q \leq p+2}\lambda_q^{2s}\|m_q\|_2 \\
\lesssim & c_r \kappa \sum_{p > Q_{b,h}}\lambda_p^{s+1} \|w_p\|_2 \sum_{q \leq p+2}\lambda_q^{s+1}\|m_q\|_2\lambda_{q-p}^{s}\lambda_q^{-1}\\
\lesssim & c_r \kappa \sum_{q \geq -1 }(\lambda_q^{2s+2}\|w_q\|_2^2+\lambda_q^{2s+2}\|m_q\|_2^2).
\end{split}
\end{equation}

\subsection{Estimation of H}
By Bony's paraproduct decomposition, one has
\begin{equation}\notag
\begin{split}
H=& \sum_{q \geq -1}\sum_{|p-q| \leq 2} \lambda_q^{2s}\int_{\T^3}  \Delta_q(m_{p}\cdot \nabla v_{\leq p-2})m_q \mathrm{d}x\\
&+\sum_{q \geq -1}\sum_{|p-q| \leq 2} \lambda_q^{2s}\int_{\T^3}  \Delta_q(m_{\leq p-2}\cdot \nabla v_p)m_q \mathrm{d}x\\
&+\sum_{q \geq -1}\sum_{p \geq q-2} \lambda_q^{2s}\int_{\T^3}  \Delta_q(\tilde m_p \cdot \nabla v_p)m_q \mathrm{d}x\\
=:& H_1+H_2+H_3.
\end{split}
\end{equation}

By the wavenumber $Q_{u,v},$ the term $H_1$ can be split into three parts.
\begin{equation}\notag
\begin{split}
H_{1}= & \sum_{-1 \leq p \leq Q_{u,v}+2}\sum_{|q-p| \leq 2} \lambda_q^{2s}\int_{\T^3}  \Delta_q(m_{p}\cdot \nabla v_{\leq p-2})m_q \mathrm{d}x\\
&+ \sum_{p > Q_{u,v}+2}\sum_{|q-p| \leq 2} \lambda_q^{2s}\int_{\T^3}  \Delta_q(m_{p}\cdot \nabla v_{\leq Q_{u,v}})m_q \mathrm{d}x\\
&+ \sum_{p > Q_{u,v}+2}\sum_{|q-p| \leq 2} \lambda_q^{2s}\int_{\T^3}  \Delta_q(m_{p}\cdot \nabla v_{(Q_{u,v}, p-2]})m_q \mathrm{d}x\\
=:& H_{11}+H_{12}+H_{13}.
\end{split}
\end{equation}

One can estimate $H_{11}$ with the help of Definition \ref{uvnbr}, H\"older's and Young's inequalities.
\begin{equation}\notag
\begin{split}
|H_{11}| \leq & \sum_{p \geq -1} \|m_p\|_\frac{2r}{r-2}\sum_{|q-p| \leq 2} \lambda_q^{2s}\|m_q\|_2 \sum_{-1 \leq p' \leq p-2 }\lambda_{p'}\| v_{p'}\|_r  \\
\lesssim & c_r \kappa \sum_{p \geq -1} \lambda_p^{\frac{3}{r}}\|m_p\|_2\sum_{|q-p| \leq 2} \lambda_q^{2s}\|m_q\|_2\sum_{-1 \leq p' \leq p-2 }\lambda_{p'}^{2-\frac{3}{r}}\\
\lesssim & c_r \kappa \sum_{p \geq -1} \lambda_p^{s+1}\|m_p\|_2\sum_{|q-p| \leq 2} \lambda_q^{s+1}\|m_q\|_2\sum_{-1 \leq p' \leq p-2 }\lambda_{p'-p}^{2-\frac{3}{r}}\\
\lesssim & c_r \kappa \sum_{q \geq -1}\lambda_q^{2s+2}\|m_q\|_2.
\end{split}
\end{equation}

To estimate $H_{12},$ one recalls Definition \ref{uvnbr} and applies H\"older's and Young's inequalities.
\begin{equation}\notag
\begin{split}
|H_{12}| \leq & \|\nabla v_{\leq Q_{u,v}}\|_r \sum_{p > Q_{u,v}} \|m_p\|_\frac{2r}{r-2}\sum_{|q-p| \leq 2} \lambda_q^{2s}\|m_q\|_2\\
\lesssim & \Lambda_{u,v}^{-1+\frac{3}{r}}\|v_{\leq Q_{u,v}}\|_r \sum_{p > Q_{u,v}} \lambda_p^{2-\frac{3}{r}}\|m_p\|_{\frac{2r}{r-2}}\sum_{|q-p| \leq 2} \lambda_q^{2s}\|m_q\|_2\\
\lesssim & c_r \kappa \sum_{p \geq -1} \lambda_p^{s+1}\|m_p\|_2\sum_{|q-p| \leq 2} \lambda_q^{s+1}\|m_q\|_2\\
\lesssim & c_r \kappa \sum_{q \geq -1}\lambda_q^{2s+2}\|m_q\|_2.
\end{split}
\end{equation}

As a result of Definition \ref{uvnbr}, H\"older's, Young's and Jensen's inequalities, one has
\begin{equation}\notag
\begin{split}
|H_{13}| \leq & \sum_{p > Q_{u,v}+2} \|m_p\|_\frac{2r}{r-2}\sum_{|q-p| \leq 2} \lambda_q^{2s}\|m_q\|_2 \sum_{Q_{u,v} < p' \leq p-2}\lambda_{p'}\|v\|_r\\
\leq & c_r \kappa \sum_{p > Q_{u,v}+2} \lambda_p^2\|m_p\|_2\sum_{|q-p| \leq 2} \lambda_q^{2s}\|m_q\|_2 \sum_{Q_{u,v} < p' \leq p-2}\lambda_{p'-p}^{2-\frac{3}{r}}\\
\lesssim & c_r \kappa \sum_{p \geq -1} \lambda_p^{s+1}\|m_p\|_2\sum_{|q-p| \leq 2} \lambda_q^{s+1}\|m_q\|_2\\
\lesssim & c_r \kappa \sum_{q \geq -1}\lambda_q^{2s+2}\|m_q\|_2.
\end{split}
\end{equation}

$H_2$ is split into lower and higher modes.
\begin{equation}\notag
\begin{split}
H_{2}= & \sum_{-1 \leq p \leq Q_{u,v}}\sum_{|q-p| \leq 2} \lambda_q^{2s}\int_{\T^3}  \Delta_q(m_{\leq p-2}\cdot \nabla v_p)m_q \mathrm{d}x\\
&+ \sum_{p > Q_{u,v}}\sum_{|q-p| \leq 2} \lambda_q^{2s}\int_{\T^3}  \Delta_q(m_{\leq p-2}\cdot \nabla v_p)m_q \mathrm{d}x\\
=:& H_{21}+H_{22},
\end{split}
\end{equation}
which are estimated by Definition \ref{uvnbr}, H\"older's, Young's and Jensen's inequalities.
\begin{equation}\notag
\begin{split}
|H_{21}| \leq & \sum_{-1 \leq p \leq Q_{u,v}}\|\nabla v_p\|_r \sum_{|q-p| \leq 2} \lambda_q^{2s}\|m_q\|_2 \sum_{p' \leq p-2}\|m_{p'}\|_\frac{2r}{r-2}\\
\leq & c_r \kappa \sum_{-1 \leq p \leq Q_{u,v}}\lambda_p^{2-\frac{3}{r}} \sum_{|q-p| \leq 2} \lambda_q^{2s}\|m_q\|_2 \sum_{p' \leq p-2}\lambda_{p'}^\frac{3}{r}\|m_{p'}\|_2\\
\leq & c_r \kappa \sum_{-1 \leq p \leq Q_{u,v}}\sum_{|q-p| \leq 2} \lambda_q^{s+1}\|m_q\|_2 \sum_{p' \leq p-2}\lambda_{p'}^{s+1}\|m_{p'}\|_2\lambda_{p'-q}^{\frac{3}{r}-s-1} \\
\lesssim & c_r \kappa \sum_{q \geq -1}\lambda_q^{2s+2}\|m_q\|_2;
\end{split}
\end{equation}
\begin{equation}\notag
\begin{split}
|H_{22}| \leq & \sum_{p > Q_{u,v}}\|\nabla v_p\|_r \sum_{|q-p| \leq 2} \lambda_q^{2s}\|m_q\|_2 \sum_{p' \leq p-2}\|m_{p'}\|_\frac{2r}{r-2}\\
\lesssim & \sum_{p > Q_{u,v}}\lambda_p^{-1+\frac{3}{r}}\|v_p\|_r \sum_{|q-p| \leq 2} \lambda_q^{2s+2-\frac{3}{r}}\|m_q\|_2 \sum_{p' \leq p-2}\lambda_{p'}^\frac{3}{r}\|m_{p'}\|_2\\
\lesssim & c_r \kappa \sum_{p > Q_{u,v}}\sum_{|q-p| \leq 2} \lambda_q^{s+1}\|m_q\|_2 \sum_{p' \leq p-2}\lambda_{p'}^{s+1}\|m_{p'}\|_2\lambda_{q-p'}^{s+1-\frac{3}{r}}\\
\lesssim & c_r \kappa \sum_{q \geq -1}\lambda_q^{2s+2}\|m_q\|_2.
\end{split}
\end{equation}

One also divides $H_3$ into two terms.
\begin{equation}\notag
\begin{split}
H_{3}= & \sum_{-1 \leq p \leq Q_{u,v}}\sum_{q \leq p+2} \lambda_q^{2s}\int_{\T^3}  \Delta_q(\tilde m_p \cdot \nabla  v_p)m_q \mathrm{d}x\\
&+ \sum_{p > Q_{u,v}}\sum_{q \leq p+2} \lambda_q^{2s}\int_{\T^3}  \Delta_q(\tilde m_p \cdot \nabla v_p)m_q \mathrm{d}x\\
=:& H_{31}+H_{32}.
\end{split}
\end{equation}

By Definition \ref{uvnbr}, H\"older's, Young's and Jensen's inequalities, one has
\begin{equation}\notag
\begin{split}
|H_{31}| \leq & \sum_{-1 \leq p \leq Q_{u,v}}\|\nabla v_p\|_r\|\tilde m_p\|_2\sum_{q \leq p+2} \lambda_q^{2s}\|m_q\|_\frac{2r}{r-2} \\
\lesssim & c_r \kappa \sum_{-1 \leq p \leq Q_{u,v}}\lambda_{p}^{2-\frac{3}{r}}\|m_p\|_2\sum_{q \leq p+2} \lambda_q^{2s+\frac{3}{r}}\|m_q\|_2\\
\lesssim & c_r \kappa \sum_{-1 \leq p \leq Q_{u,v}}\lambda_{p}^{s+1}\|m_p\|_2\sum_{q \leq p+2} \lambda_q^{s+1}\|m_q\|_2\lambda_{q-p}^{s-1+\frac{3}{r}}\\
\lesssim & c_r \kappa \sum_{q \geq -1}\lambda_q^{2s+2}\|m_q\|_2.
\end{split}
\end{equation}

For $H_{32},$ the following estimate holds.
\begin{equation}\notag
\begin{split}
|H_{32}| \leq & \sum_{p >Q_{u,v}}\|\nabla v_p\|_r\|\tilde m_p\|_2\sum_{q \leq p+2} \lambda_q^{2s}\|m_q\|_\frac{2r}{r-2} \\
\lesssim & c_r \kappa \sum_{p > Q_{u,v}}\lambda_{p}^{2-\frac{3}{r}}\|m_p\|_2\sum_{q \leq p+2} \lambda_q^{2s+\frac{3}{r}}\|m_q\|_2\\
\lesssim & c_r \kappa \sum_{p > Q_{u,v}}\lambda_{p}^{s+1}\|m_p\|_2\sum_{q \leq p+2} \lambda_q^{s+1}\|m_q\|_2\lambda_{q-p}^{s-1+\frac{3}{r}}\\
\lesssim & c_r \kappa \sum_{q \geq -1}\lambda_q^{2s+2}\|m_q\|_2.
\end{split}
\end{equation}

\subsection{Conclusion}
As the terms $I$ and $J$ are already estimated in Section \ref{EMH}, one sums up all the previous estimates and chooses a suitable constant $c_r$ to obtain
\begin{equation}\notag
\frac{\mathrm{d}}{\mathrm{d}t}\sum_{q \geq -1}\big(\|w_q\|_2^2+\|m_q\|_2^2\big)\lesssim -\sum_{q \geq -1}\lambda_q^{2}\big(\|w_q\|_2^2+\|m_q\|_2^2\big)\lesssim \sum_{q \geq -1}\big(\|w_q\|_2^2+\|m_q\|_2^2\big).
\end{equation}
One can see that $\big(\|w\|_{L^2}^2+\|m\|_{L^2}^2\big)$ decays to $0$ exponentially as $t \to \infty$ as a result of Gr\"onwall's inequality.

\cbdu

\section{Bounds on the wavenumbers}

In \cite{CDK}, it was shown that the time average of the determining wavenumber for a weak solution to the Navier-Stokes equations is bounded above by Kolmogorov's dissipation wavenumber via the average energy dissipation rate $\varepsilon:= \langle \|\nabla u\|^2_{L^2} \rangle,$ where $\langle \cdot \rangle$ signifies the time average. For the 2D MHD system, it is also known that explicit dimension estimates of functional invariant sets can be given by the energy dissipation rate.

Yet, in the case of the Hall-MHD system, it seems impossible to bound the wavenumber $\Lambda_{b,h}(t)$ using the average magnetic energy dissipation rate $\langle \|\nabla b\|_{L^2}^2  \rangle$. Fortunately, restricting one's attentions to strong solutions can lead to a reasonable bound on $\Lambda_{b,h}(t)$ in an average sense. Indeed, whenever $\Lambda_{b,h}(t)> \lambda_0$, it must be that one of the conditions in Definition \ref{wvnbr} is unfulfilled, i.e.,  $\|b_{Q_{b,h}(t)}\|_\infty > c_r \kappa$ or $\|b_{\leq Q_{b,h}(t)-1}\|_\infty > c_r \kappa.$

The inequality $\|b_{Q_{b,h}(t)}\|_\infty > c_r \kappa$ implies that $$\Lambda_{b,h}(t)\|b_{Q_{b,h}(t)}\|_\infty > c_r \kappa \Lambda_{b,h}(t).$$
By Lemma \ref{brn}, one has $$ \|\nabla b\|_\infty^2 \geq \|\nabla b_{Q_{b,h}(t)}\|^2_\infty > \big(c_r \kappa \Lambda_{b,h}(t)\big)^2.$$
Meanwhile, if $\|b_{\leq Q_{b,h}(t)-1}\|_\infty > c_r \kappa,$ then
$$\Lambda_{b,h}(t)\|b_{\leq Q_{b,h}(t)-1}\|_\infty > c_r \kappa \Lambda_{b,h}(t),$$
which, by Lemma \ref{brn}, yields
$$ \|\nabla b\|_\infty^2 \geq \|\nabla b_{\leq Q_{b,h}(t)-1}\|^2_\infty > \big(c_r \kappa \Lambda_{b,h}(t)\big)^2.$$
Hence, for $(u,b) \in L^\infty\big(0,\infty; (H^s(\T^3))^2\big),$ one has, by Theorem \ref{rgc}, the following bound.
$$\langle \Lambda_{b,h}^2 \rangle \lesssim \|\nabla b\|_{L^2(0,T; L^\infty(\R^3))} < \infty.$$

\end{document}